\RequirePackage{tikz}
\usetikzlibrary{positioning, quotes, decorations.pathmorphing,shapes, arrows, arrows.meta}

\documentclass[pdflatex,sn-mathphys]{sn-jnl}

\usepackage{array}
\usepackage{bm,eucal}
\usepackage{amssymb,mathrsfs, amsmath, adjustbox}

\jyear{2023}

\theoremstyle{thmstyleone}
\newtheorem{theorem}{Theorem}
\newtheorem{proposition}[theorem]{Proposition}
\newtheorem{lemma}[theorem]{Lemma}
\newtheorem{corollary}[theorem]{Corollary}

\theoremstyle{thmstyletwo}
\newtheorem{remark}{Remark}
\newtheorem{example}{Example}
\newtheorem{notation}{Notation}

\theoremstyle{thmstylethree}
\newtheorem{definition}{Definition}

\usepackage{soul}

\begin{document}

\title[Weak Compactness in $ W^{k, 1} $ with the Existence of Minimizers]{Weak Compactness Criterion in $ W^{k, 1} $ with an Existence Theorem of Minimizers}

\author[1]{\fnm{Cheng} \sur{Chen}}\email{chenchengscu@outlook.com}
\author[2]{\fnm{Mattie} \sur{Ji}}\email{mji13@sas.upenn.edu}
\author[1]{\fnm{Yan} \sur{Tang}}\email{tangyan@ctbu.edu.cn}
\author*[1]{\fnm{Shiqing} \sur{Zhang}}\email{zhangshiqing@scu.edu.cn}

\affil*[1]{\orgdiv{College of Mathematics}, \orgname{Sichuan University}, \orgaddress{\street{South Section 1 of the 1st Ring Road}, \city{Chengdu}, \postcode{610064}, \state{Sichuan}, \country{China}}}

\affil*[2]{\orgdiv{Department of Mathematics}, \orgname{University of Pennsylvania}, \orgaddress{\street{209 South 33rd St}, \city{Philadelphia}, \postcode{19104}, \state{Pennsylvania}, \country{USA}}}

\abstract{There is a rich theory of existence theorems for minimizers over reflexive Sobolev spaces (ex. Eberlein-\v{S}mulian theorem). However, the existence theorems for many variational problems over non-reflexive Sobolev spaces remain underexplored. In this paper, we investigate various examples of functionals over non-reflexive Sobolev spaces. To do this, we prove a weak compactness criterion in $W^{k,1}$ that generalizes the Dunford-Pettis theorem, which asserts that relatively weakly compact subsets of $ L^1 $ coincide with equi-integrable families. As a corollary, we also extend an existence theorem of minimizers from reflexive Sobolev spaces to non-reflexive ones. This work is also benefited and streamlined by various concepts in category theory. 
}
% \abstract{Nelson Dunford and Billy James Pettis [{\em Trans. Amer. Math. Soc.}, 47 (1940), pp. 323--392] proved that relatively weakly compact subsets of $ L^1 $ coincide with equi-integrable families. We expand it to the case of $ W^{k,1} $ - the non-reflexive Sobolev spaces - by a tailor-made isometric operator. Herein we extend an existence theorem of minimizers from reflexive Sobolev spaces to non-reflexive ones.}

\keywords{Dunford-Pettis theorem, non-reflexive Sobolev space, weak compactness in $ W^{k,1} $, existence theorem of a minimizer, calculus of variations.}

\pacs[MSC Classification]{46N20, 46E35, 28A20, 46E30}

\maketitle

\section{Introduction and Main Results}\label{sec1}

Historically, starting from Leonhard Euler and Joseph-Louis Lagrange, a series of practical problems in natural science including Fermat's principle, the brachistochrone problem, and Dirichlet's principle \cite[see][]{ZhangShiqing,Brezis,Courant,Courant&Hilbert}, can be transformed into minimizing a functional of the following form:
\[
	J(u) = \int_{\Omega} f\bigl(x, u(x), \nabla u(x)\bigr) \mathrm{d}x,
\]
where $ \Omega $ is a measurable space that satisfies some smoothness conditions, and $ u $ varies in a suitable function space, $ f $ is uniformly integrable on $ \Omega $ and known as the Lagrangian state function.

The methods for solving the above problems are collectively called the \textit{calculus of variations}. These minimizers were initially considered to exist ``naturally'' until Karl Weierstrass constructed a counterexample \cite[see][Example 4.6 on page 122]{Dacorogna} in 1870. Since then, people have studied and obtained many existence theorems of minimizers. For example, Dirichlet's principle was first proved in 1899 (with some strong conditions) then later extended to more general cases by David Hilbert \cite{Hilbert} in 1904.

To be more precise, a real-valued functional $ J $ induces a preorder $ \preceq $ over its domain such that for all $ u, v $ in a Banach space $ E := \mathrm{dom}(J) $: 
\[
	u \preceq v \text{ if and only if } J(u) \geq J(v).
\]
A net $ \{u_\alpha\}_{\alpha \in \Lambda} $ equipped with the preorder above is called a decreasing net of $ J $. And whether a minimizer exists is equivalent to whether there is a greatest element for each of these decreasing nets. Therefore we should pay attention to whether any decreasing sequence $ \{u_n\} $ satisfies that
\[
	\lim_{n \to + \infty} J(u_n) = \inf_{u \in E} J(u)
\]
(which is so called the minimizing sequence of $ J $) can reach the infimum. Considering that the real line has a natural topology and selecting ``what are the continuous functions from $ E $ to $ \mathbb{R} $'', the topology on $ \mathbb{R} $ can be induced to $ E $ (e.g. the norm topology and the weak topology) while keeping $ J $ lower semicontinuous under this topology. Since lower semicontinuous functions can reach the infimum on compact sets, if a minimizing sequence $ \{ u_n \} $ is contained in some compact set, then there exists some $ \bar{u} $ in the closure of $ \{ u_n \} $ such that the functional $ J $ can reach the infimum at $ \bar{u} $.

In fact, the existence of minimizers depends on what domain space one chooses. Roughly speaking, the ``larger'' the space, the weaker the topology, the ``fewer'' open sets, and the easier it is for a set to be relatively compact. For aesthetic and practical reasons, one may typically expect a variational problem to have a solution as smooth as possible. However, it is usually difficult to directly find a solution with high regularity in a space with a strong topological structure. But if (weak) solutions can exist in a ``larger'' space, one only needs to verify that they are regular. For regularity, people have also developed many profound theories \cite[see][Section 8.3 on page 458]{Evans} and \cite[Section 7.3 on page 191 and Section 9.6 on page 298]{Brezis}, but we will only focus on existence. Arguably, Sobolev spaces \cite[see][3.2 on page 59]{Adams} are tailor-made for this approach. By replacing classical derivatives with weak generalized derivatives, Sobolev spaces are ``large'' enough. But they are not too ``bad'' since even the smooth function space $ C^{\infty} $ can be densely embedded in them (the Meyers-Serrin theorem \cite{Meyers}).

Assuming a functional is (weakly) lower semicontinuous in some Sobolev space, one would expect its corresponding minimizing sequence to be (weakly) relatively compact. For reflexive Sobolev spaces, we can apply the Eberlein-\v{S}mulian theorem \cite[see][page 141]{Yosida}, that is, a Banach space is reflexive if and only if any bounded set inside is weakly relatively compact. Based on this and James's theorem \cite{James}, Tang, Zhang, and Guo \cite{Tangyan} have made a series of discussions on the relationship between the reflexivity of a Banach space and the existence of minimizers of a sequentially weakly lower semicontinuous functional on it. Using the Eberlein-\v{S}mulian theorem, researchers have obtained some existence theorems of functional minimizers on reflexive Sobolev spaces. An example is Theorem 3.30 in \cite[page 106]{Dacorogna}, which can be used to prove Dirichlet's principle as an application.

However, there are still various important practical problems whose integral functionals are defined on non-reflexive Sobolev spaces, such as the Plateau problem \cite[see][]{Courant,CourantR} (which is to show the existence of a minimal surface with a given boundary). Another simpler example follows: 

\begin{example}\cite[Example 4.5 on page 122 and Example 4.9 on page 124]{Dacorogna}\label{example1}
	Let $ u $ vary in the following function space
	\[
		E = \{u \in W^{1, 1}(0,1) \vert \: u(0) = 0, \: u(1) = 1\},
	\]
	where the value of $ u $ on the boundary is defined in the sense of the trace operator. Let $J_1$ and $J_2$ be the functionals
	\[
		J_1(u) = \int_{0}^{1} \sqrt{u^2(t) + \dot{u}^2(t)} \mathrm{d}t
	\]
	and
	\[
		J_2(u) = \int_{0}^{1} \lvert \dot{u}(t) \rvert \mathrm{d}t.
	\]
	Can their respective infimum be attained in $E$?
\end{example}

Since the Eberlein-\v{S}mulian theorem no longer applies in this case, we need something stronger than boundedness to regain relative weak compactness. For the specific function space $ L^1 $, Nelson Dunford and Billy James Pettis obtained the following theorem in 1940.

\begin{theorem}[Dunford-Pettis theorem: weak compactness criterion in $ L^1 $]\textup{\cite[Theorem 3.2.1 on page 376]{Dunford2}, \cite[see also][Problem 23 on page 466]{Brezis}}\label{Dunford-Pettis}
	Suppose that
	\begin{description}
		\item[$ \mathbf{H_{\boldsymbol{\sigma}}^1}{[\Omega, \mathcal{F}]} $ \textup{(hypothesis):}] 
		$ \Omega $ is a $ \sigma $-finite measure space, and $ \mathcal{F} $ is a subset of $ L^1(\Omega) $.
	\end{description}

	Then the following statements $ \mathbf{B^1} \wedge \mathbf{EI^1} \wedge \mathbf{E_{\infty}^1} $ and $ \mathbf{C_w^1} $ are equivalent. Here, the two statements have the following meaning.
 
	\begin{description}
		\item[1. $ \mathbf{B^1} \wedge \mathbf{EI^1} \wedge \mathbf{E_{\infty}^1}{[\Omega, \mathcal{F}, \exists M, \exists \delta, \exists \omega]} $\textup{:}] 
		$ \mathcal{F} $ is an equi-integrable family.
		\item[2. $ \mathbf{C_w^1}{[\Omega, \mathcal{F}]} $ \textup{(weak compactness in $ L^1 $):}] 
		the set $ \mathcal{F} $ is relatively weakly compact in $ L^1(\Omega) $.
	\end{description}
\end{theorem}

By an ``equi-integrable family" in Theorem~\ref{Dunford-Pettis}, we mean the following.

\begin{definition}[equi-integrable family]\textup{\cite[4.36 on page 129]{Brezis}}
	A subset $ \mathcal{F} \subset L^1(\Omega) $ is said to be equi-integrable if it satisfies the following three conditions:
 
	\begin{description}
		\item[a. $ \mathbf{B^1}{[\Omega, \mathcal{F}, \exists M]} $ \textup{(boundedness in $ L^1 $):}] 
		there exists $ M \in \mathbb{R}_{> 0} $ such that $ \lVert f \rVert_{L^1(\Omega)} \leq M $ for every $ f \in \mathcal{F} $.
  
		\item[b. $ \mathbf{EI^1}{[\Omega, \mathcal{F}, \exists \delta]} $ \textup{(equi-integrability):}] 
		there exists a function $ \delta \colon \mathbb{R}_{> 0} \to \mathbb{R}_{> 0} $ such that
		\[
			\int_A \lvert f \rvert \mathrm{d}\lambda < \varepsilon
		\]
		for each $ f \in \mathcal{F} $, for any $ \varepsilon \in \mathbb{R}_{> 0} $, and for all $ A $ which is measurable in $ \Omega $ with its measure $ \lambda(A) < \delta(\varepsilon) $.
  
		\item[c. $ \mathbf{E_{\infty}^1}{[\Omega, \mathcal{F}, \exists \omega]} $ \textup{(equi-integrability at infinity):}] 
		there exists a set-valued mapping $ \omega $ from $ \mathbb{R}_{> 0} $ to $ 2^{\Omega} $ (the power set of $ \Omega $) such that $ \omega(\varepsilon) $ is measurable with its measure $ \lambda \circ \omega(\varepsilon) $ being finite, and there satisfies that
		\[
			\int_{\Omega \setminus \omega(\varepsilon)} \lvert f \rvert \mathrm{d}\lambda < \varepsilon
		\]
		for every $ f \in \mathcal{F} $ and any $ \varepsilon \in \mathbb{R}_{> 0} $.
	\end{description}
	The notation $ \lambda $ here and thereafter refers to the Lebesgue measures \cite[page 14]{Adams} in real Euclidean spaces. Note that when $ \Omega $ is bounded, the third condition $ \mathbf{E_{\infty}^1}{[\Omega, \mathcal{F}, \exists \omega]}$ is naturally true by choosing $ \omega \equiv \Omega $.
\end{definition}

\begin{remark}
	We mark the hypothesis and statements in the Dunford-Pettis theorem with symbols for two reasons:
    \begin{description}
        \item[1.] We do not have to repeat the sentence from beginning to end every time. For example, we can shrink the Dunford-Pettis theorem to the following form:
	\[
		\mathbf{H_{\boldsymbol{\sigma}}^1} \wedge \mathbf{B^1} \wedge \mathbf{EI^1} \wedge \mathbf{E_{\infty}^1}[\Omega, \mathcal{F}, \exists M, \exists \delta, \exists \omega] \Leftrightarrow \mathbf{H_{\boldsymbol{\sigma}}^1} \wedge \mathbf{C_w^1}[\Omega, \mathcal{F}].
	\]
	We also use superscripts and subscripts to express some key information. For the full list of the types of symbols that will show up in this paper and their meaning, please see Appendix~\ref{notations}.
        \item[2.] We can put the elements necessary for a sentence to have a specific meaning into the following brackets and emphasize them. For example, in $ \mathbf{B^1}[\Omega, \mathcal{F}, \exists M] $, we know that if we fix $ \Omega $ and $ \mathcal{F} $ and can prove the existence of $ M $, then $ \mathbf{B^1}[\Omega, \mathcal{F}, \exists M] $ is a certain sentence. In fact, we can construct a thin category whose objects are sentences and morphisms are derivation symbols ``$ \Rightarrow $''. It does not matter if we directly treat ``$ \Rightarrow $'' as ``deduces'' and ``$ \Leftrightarrow $'' as ``is equivalent to''. Obviously ``$ \Leftrightarrow $'' are the isomorphisms in this category. See Appendix \ref{Details} for more details. 
    \end{description}
\end{remark}

Based on the Dunford-Pettis theorem above, we will prove the following theorem on the non-reflexive Sobolev space $W^{k,1}$ through a technique of ``disassembly and assembly'':

\begin{theorem}[weak compactness criterion in $ W^{k, 1} $]\label{T1}
	Suppose that
	\begin{description}
		\item[$ \mathbf{H^{k, 1}_{n.o.}}{[\Omega, \mathcal{F}]} $ \textup{(hypothesis):}] 
		$ \Omega $ is a non-empty open set in $ \mathbb{R}^d $, and $ \mathcal{F} $ is a subset of $ W^{k,1}(\Omega) $.
	\end{description}
	
	Then the two statements $ \mathbf{B^{k, 1}} \wedge \mathbf{EI^{k, 1}} \wedge \mathbf{E_{\infty}^{k, 1}} $ and $ \mathbf{C_w^{k, 1}} $ are equivalent. Here, the two statements have the following meaning.
 
	\begin{description}
		\item[1. $ \mathbf{B^{k, 1}} \wedge \mathbf{EI^{k, 1}} \wedge \mathbf{E_{\infty}^{k, 1}}{[\Omega, \mathcal{F}, \exists M, \exists \delta, \exists \omega]} $\textup{:}] 
		$ \mathcal{F} $ is said to be equi-integrable under the sense of $ W^{k,1} $ provided the following three conditions:

		\begin{description}
			\item[a. $ \mathbf{B^{k, 1}}{[\Omega, \mathcal{F}, \exists M]} $ \textup{(boundedness in $ W^{k, 1} $):}] 
			there exists a constant $ M \in \mathbb{R}_{> 0} $ such that $ \lVert f \rVert_{W^{k,1}(\Omega)} \leq M $ for every $ f \in \mathcal{F} $.
			\item[b. $ \mathbf{EI^{k, 1}}{[\Omega, \mathcal{F}, \exists \delta]} $ \textup{(equi-integrability):}] 
			there exists a function $ \delta \colon \mathbb{R}_{> 0} \to \mathbb{R}_{> 0} $ such that
			\[
				\sum_{\lvert \alpha \rvert \leq k} \int_A \lvert D^{\alpha} f \rvert \mathrm{d}\lambda < \varepsilon
			\]
			for every $ f \in \mathcal{F} $, for any $ \varepsilon \in \mathbb{R}_{> 0} $, and for all $ A $ which is measurable in $ \Omega $ with its measure $ \lambda(A) < \delta(\varepsilon) $. Here and henceforth, $\alpha$ denotes a multi-index (Notation~\ref{multi-indices}).
		\end{description}

		\begin{description}
			\item[c. $ \mathbf{E_{\infty}^{k, 1}}{[\Omega, \mathcal{F}, \exists \omega]} $ \textup{(equi-integrability at infinity):}] 
			there exists a set-valued mapping $ \omega \colon \mathbb{R}_{> 0} \to 2^{\Omega} $ such that $ \omega(\varepsilon) $ is measurable with its measure $ \lambda \circ \omega(\varepsilon) $ being finite, and there satisfies that
			\[
				\sum_{\lvert \alpha \rvert \leq k} \int_{\Omega \setminus \omega(\varepsilon)} \lvert D^{\alpha} f \rvert \mathrm{d}\lambda < \varepsilon
			\]
			for every $ f \in \mathcal{F} $ and any $ \varepsilon \in \mathbb{R}_{> 0} $.
		\end{description}
		\item[2. $ \mathbf{C_w^{k, 1}}{[\Omega, \mathcal{F}]} $ \textup{(weak compactness in $ W^{k, 1} $):}] 
		$ \mathcal{F} $ is relatively weakly compact in $ W^{k,1}(\Omega) $.
	\end{description}
\end{theorem}

When $ k = 0 $, Theorem \ref{T1} is exactly the Dunford-Pettis theorem if one regards the notation $ W^{0,1} $ as $ L^1 $.

When $ \Omega $ is bounded, the statement $ \mathbf{E_{\infty}^{k, 1}} $ in Theorem \ref{T1} is always satisfied by taking $ \omega(\varepsilon) \equiv \Omega $. Moreover, if $ \Omega $ satisfies the cone condition \cite[see][Paragraph 4.6 on page 82]{Adams}, we have the following theorem, which can also be seen as a high-dimensional generalization of Proposition 3.1.4 in \cite{Fathi}.

\begin{theorem}\label{T1'}
    Suppose that
	\begin{description}
		\item[$ \mathbf{H^{k, 1}_{b.o.cone}}{[\Omega, \mathcal{F}]} $ \textup{(hypothesis):}] 
		$ \Omega $ is a bounded open set in $ \mathbb{R}^d $, and $ \Omega $ satisfies the cone condition. $ \mathcal{F} $ is a subset of $ W^{k,1}(\Omega) $.
	\end{description}
	
	Then the statement $ \mathbf{B^{k, 1}} $ in Theorem \textup{\ref{T1}} is equivalent to the following statement $ \mathbf{B^{\lvert \alpha \rvert \in \{0, k\},1}} $.
    
    \begin{description}
			\item[$ \mathbf{B^{\lvert \alpha \rvert \in \{0, k\},1}}{[\Omega, \mathcal{F}, \exists M]} $ \textup{:}] 
			there exists a constant $ M \in \mathbb{R}_{> 0} $ such that 
                \[
				\lVert f \rVert_{L^1(\Omega)} + \sum_{\lvert \alpha \rvert = k} \lVert D^{\alpha} f \rVert_{L^1(\Omega)} \leq M
			\]
                for every $ f \in \mathcal{F} $.
    \end{description}

    The statement $ \mathbf{EI^{k, 1}} $ in Theorem \textup{\ref{T1}} can be replaced by the following statement $ \mathbf{EI^{\lvert \alpha \rvert = k,1}} $.
    \begin{description}
			\item[$ \mathbf{EI^{\lvert \alpha \rvert = k,1}}{[\Omega, \mathcal{F}, \exists \delta]} $ \textup{:}] 
			there exists a function $ \delta \colon \mathbb{R}_{> 0} \to \mathbb{R}_{> 0} $ such that
			\[
				\sum_{\lvert \alpha \rvert = k} \int_A \lvert D^{\alpha} f \rvert \mathrm{d}\lambda < \varepsilon
			\]
			for every $ f \in \mathcal{F} $, for any $ \varepsilon \in \mathbb{R}_{> 0} $, and for all $ A $ which is measurable in $ \Omega $ with its measure $ \lambda(A) < \delta(\varepsilon) $.
		\end{description}
\end{theorem}

Excitingly, this ``generalized model'' can also be ``translated'' to other well-known theorems. For example, we will generalize the Kolmogorov-M.Riesz-Fr\'echet theorem \cite[see][Corollary 4.27 on page 113]{Brezis} to obtain the following Corollary \ref{cor1} and Corollary \ref{cor1'} as follows. After this, we will generalize the Ascoli-Arzel\`a theorem \cite[see][Theorem 4.25 on page 111]{Brezis} into Corollary \ref{cor2} and Corollary \ref{cor2'}. See Appendix \ref{Poc} for the proofs of these corollaries.

\begin{corollary}[compactness criterion in $ W^{k,p} $]\label{cor1}
	Suppose that
	\begin{description}
		\item[$ \mathbf{H^{k, p}}{[\mathbb{R}^d, \mathcal{F}]} $ \textup{(hypothesis):}] 
		$ \mathcal{F} $ is a subset of $ W^{k,p}(\mathbb{R}^d) $ with $ 1 \leq p < + \infty $.
	\end{description}
	
	Then the statements $ \mathbf{B^{k, p}} \wedge \mathbf{EI_{\boldsymbol{\tau}}^{k, p}} \wedge \mathbf{E_{\infty}^{k, p}} $ and $ \mathbf{C^{k, p}} $ are equivalent. Here, the two statements have the following meaning.
 
	\begin{description}
		\item[1. $ \mathbf{B^{k, p}} \wedge \mathbf{EI_{\boldsymbol{\tau}}^{k, p}} \wedge \mathbf{E_{\infty}^{k, p}}{[\mathbb{R}^d, \mathcal{F}, \exists M, \exists \delta, \exists \omega]} $\textup{:}] 
		the following three are satisfied: 
		\begin{description}
			\item[a. $ \mathbf{B^{k, p}}{[\mathbb{R}^d, \mathcal{F}, \exists M]} $ \textup{(boundedness in $ W^{k, p} $):}] 
			there exists a constant $ M \in \mathbb{R}_{> 0} $ such that $ \lVert f \rVert_{W^{k,p}(\mathbb{R}^d)} \leq M $ for every $ f \in \mathcal{F} $.
			\item[b. $ \mathbf{EI_{\boldsymbol{\tau}}^{k, p}}{[\mathbb{R}^d, \mathcal{F}, \exists \delta]} $\textup{:}] 
			there exists a function $ \delta \colon \mathbb{R}_{> 0} \to \mathbb{R}_{> 0} $ such that
			\[
				\lVert \tau_h f - f \rVert_{W^{k, p}(\mathbb{R}^d)} < \varepsilon
			\]
			for every $ f \in \mathcal{F} $, for any $ \varepsilon \in \mathbb{R}_{> 0} $, and for all $ h \in \mathbb{R}^d $ with its norm $ \lVert h \rVert_{\mathbb{R}^d} < \delta(\varepsilon) $. Here, the notation $\tau_h f(x) := f(x+h)$ is a shift of the original function $f$ by the vector $h$.
			\item[c. $ \mathbf{E_{\infty}^{k, p}}{[\mathbb{R}^d, \mathcal{F}, \exists \omega]} $\textup{:}] 
			there exists a set-valued mapping $ \omega \colon \mathbb{R}_{> 0} \to 2^{\mathbb{R}^d} $ such that $ \omega(\varepsilon) $ is measurable with its measure $ \lambda \circ \omega(\varepsilon) $ being finite, and there satisfies that
			\[
				(\sum_{\lvert \alpha \rvert \leq k} \int_{\mathbb{R}^d \setminus \omega(\varepsilon)} \lvert D^{\alpha} f \rvert^p \mathrm{d}\lambda)^{\frac{1}{p}} < \varepsilon
			\]
			for every $ f \in \mathcal{F} $ and any $ \varepsilon \in \mathbb{R}_{> 0} $.
		\end{description}
		\item[2. $ \mathbf{C^{k, p}}{[\mathbb{R}^d, \mathcal{F}]} $ \textup{(compactness in $ W^{k, p} $):}] 
		the set $ \mathcal{F} $ is precompact in $ W^{k, p}(\mathbb{R}^d) $.
	\end{description}
\end{corollary}

\begin{corollary}\label{cor1'}
    Suppose that
	\begin{description}
		\item[$ \mathbf{H^{k, p}_{b,o,0}}{[\Omega, \mathcal{F}]} $ \textup{:}] 
		$ \Omega $ is a bounded open set in $ \mathbb{R}^n $, and $ \mathcal{F} $ is a subset of $ W_0^{k,p}(\Omega) $ with $ 1 \leq p < + \infty $.
	\end{description}
	
	Then the statement $ \mathbf{B^{k, p}} $ in Corollary \textup{\ref{cor1}} is equivalent to $ \mathbf{B^{\lvert \alpha \rvert \in \{0, k\},p}} $ as follows.
    
    \begin{description}
			\item[$ \mathbf{B^{\lvert \alpha \rvert \in \{0, k\},p}}{[\Omega, \mathcal{F}, \exists M]} $ \textup{:}] 
			there exists a constant $ M \in \mathbb{R}_{> 0} $ such that 
                \[
				\lVert f \rVert_{L^p(\Omega)} + \sum_{\lvert \alpha \rvert = k} \lVert D^{\alpha} f \rVert_{L^p(\Omega)} \leq M
			\]
                for every $ f \in \mathcal{F} $.
    \end{description}

    The statement $ \mathbf{EI_{\boldsymbol{\tau}}^{k, p}} $ in Corollary \textup{\ref{cor1}} can be replaced by the following $ \mathbf{EI_{\boldsymbol{\tau}}^{\lvert \alpha \rvert = k,p}} $.
    \begin{description}
			\item[$ \mathbf{EI_{\boldsymbol{\tau}}^{\lvert \alpha \rvert = k,p}}{[\Omega, \mathcal{F}, \exists \delta]} $\textup{:}] 
			there exists a function $ \delta \colon \mathbb{R}_{> 0} \to \mathbb{R}_{> 0} $ such that
			\[
				\sum_{\lvert \alpha \rvert = k} \lVert \tau_h (D^{\alpha}f) - D^{\alpha}f \rVert_{L^p(\mathbb{R}^d)} < \varepsilon
			\]
			for every $ f \in \mathcal{F} $, for any $ \varepsilon \in \mathbb{R}_{> 0} $, and for all $ h \in \mathbb{R}^d $ with its norm $ \lVert h \rVert_{\mathbb{R}^d} < \delta(\varepsilon) $. Here, the functions $ D^{\alpha}f $ are extended to be $0$ outside $\Omega$.
		\end{description}
\end{corollary}

\begin{corollary}[compactness criterion in $ C^m $]\label{cor2}
	Suppose that
	\begin{description}
		\item[$ \mathbf{H^m}{[K, \mathcal{F}]} $ \textup{(hypothesis):}] 
		$K$ is a compact metric space which satisfies the uniform $ C^m $ regularity condition \textup{\cite[see][4.10 on page 84]{Adams}}, and $ \mathcal{F} $ is a subset of $ C^m(K) $.
	\end{description}
	
	Then the two statements $ \mathbf{B^m} \wedge \mathbf{EC^m} $ and $ \mathbf{C^m} $ are equivalent. Here, the two statements have the following meaning.
	\begin{description}
		\item[1. $ \mathbf{B^m} \wedge \mathbf{EC^m}{[K, \mathcal{F}, \exists M, \exists \delta]} $\textup{:}] 
		$ \mathcal{F} $ is said to be uniformly equicontinuous under the sense of $ C^m $ provided the following two conditions:

		\begin{description}
			\item[a. $ \mathbf{B^m}{[K, \mathcal{F}, \exists M]} $ \textup{(boundedness in $ C^m $):}] 
			there exists a constant $ M \in \mathbb{R}_{> 0} $ such that $ \lVert f \rVert_{C^m(K)} := \max_{\lvert \alpha \rvert \leq m} \lVert D^{\alpha}f \rVert_{C(K)} \leq M $ for every $ f \in \mathcal{F} $.
			\item[b. $ \mathbf{EC^m}{[K, \mathcal{F}, \exists \delta]} $ \textup{(equicontinuity):}] 
			there exists a function $ \delta \colon \mathbb{R}_{> 0} \to \mathbb{R}_{> 0} $ such that
			\[
				\max_{\lvert \alpha \rvert \leq m} \lvert (D^{\alpha}f)(x_1) - (D^{\alpha}f)(x_2) \rvert < \varepsilon
			\]
			for every $ f \in \mathcal{F} $, for any $ \varepsilon \in \mathbb{R}_{> 0} $, and for all $ x_1 $ and $ x_2 $ in $ K $ with the distance between them $ d(x_1, x_2) < \delta(\varepsilon) $.
		\end{description}
		\item[2. $ \mathbf{C^m}{[K, \mathcal{F}]} $ \textup{(compactness in $ C^m $):}] 
		the set $ \mathcal{F} $ is precompact in $ C^m(K) $.
	\end{description}
\end{corollary}

\begin{corollary}\label{cor2'}
    Statement $ \mathbf{EC^m} $ in Corollary \textup{\ref{cor2}} can be replaced by the following $ \mathbf{EC^{\lvert \alpha \rvert = m}} $.
    \begin{description}
			\item[$ \mathbf{EC^{\lvert \alpha \rvert = m}}{[K, \mathcal{F}, \exists \delta]} $\textup{:}] 
			there exists a function $ \delta \colon \mathbb{R}_{> 0} \to \mathbb{R}_{> 0} $ such that
			\[
				\max_{\lvert \alpha \rvert = m} \lvert (D^{\alpha}f)(x_1) - (D^{\alpha}f)(x_2) \rvert < \varepsilon
			\]
			for every $ f \in \mathcal{F} $, for any $ \varepsilon \in \mathbb{R}_{> 0} $, and for all $ x_1 $ and $ x_2 $ in $ K $ with the distance between them $ d(x_1, x_2) < \delta(\varepsilon) $.
		\end{description}
\end{corollary}

\begin{remark}
    For an arbitrary $ K $ satisfying $ \mathbf{H^m} $ in Corollary \ref{cor2}, we have not been able to generalize $ \mathbf{B^m} $ to only verify the boundedness of $ D^{\alpha}\mathcal{F} $ when $ \lvert \alpha \rvert $ is $0$ or $m$. We considered that continuous functions are very similar to Lebesgue space $ L^{\infty} $, and $ L^{\infty} $ on finite measure space can be regarded as the intersection of all $ L^p $ spaces with $ p \in [1,+\infty) $. For each $ p $, there exists a constant $ M_p $ such that
    \[
    \lVert f \rVert_{L^p(K)} + \sum_{\lvert \alpha \rvert = k} \lVert D^{\alpha} f \rVert_{L^p(K)} \leq M_p
    \]
    for every $ f \in \mathcal{F} $. However $ \sup_{p \in [1,+\infty)} M_p $ might be $ + \infty $.

     When $ K $ is a subset of $ \mathbb{R} $, the Landau–Kolmogorov inequality \cite{Landau} can be considered to apply. But our topic is not interpolation inequality so we will not go into it too much.
\end{remark}

With the compactness criterion in $ W^{k, 1} $, we can extend some theorems that were previously proven in reflexive spaces to some non-reflexive spaces. As an application, we will prove the following theorem as an extension of Theorem 3.30 in \cite[page 106]{Dacorogna}: 

\begin{theorem}[an existence theorem of minimizers in $ W^{1, 1} $]\label{T2}
	Let $ \Omega $ be a bounded open subset of $ \mathbb{R}^d $, and $ \Omega $ satisfies the cone condition. Let $ L \colon \Omega \times \mathbb{R}^N \times \mathbb{R}^{N \times d} \to \mathbb{R} \cup \{+ \infty\} $ be a Carath\'{e}odory function \textup{\cite[see][Definition 3.5 on page 75]{Dacorogna}} satisfying
	\begin{equation}\label{E4.1}
		L(x,u,\xi) \geq a \lvert \xi \rvert + c(x)
	\end{equation}
	for almost every $ x \in \Omega $ and for every $ (u, \xi) \in \mathbb{R}^N \times \mathbb{R}^{N \times d} $ and for some $ a \in \mathbb{R}_{> 0} $, $ c \in L^1(\Omega) $. And let the functional
	\[
		J(u) = \int_{\Omega} L\bigl(x, u(x), \nabla u(x)\bigr) \mathrm{d}x
	\]
	be finite at some $ u_0 \in W^{1, 1}(\Omega) $. Then every minimizing sequence in $ u_0 + W_0^{1, 1}(\Omega) $ is bounded.
	
	In addition, if $ \xi \mapsto L(x,u,\xi) $ is convex, and a minimizing sequence $ \{u_n\} $ is found to satisfy the statement $ \mathbf{B^{1,1}} \wedge \mathbf{EI^{\lvert \alpha \rvert = 1}}{[\Omega, \{u_n\}, \exists M, \exists \delta]} $. Then $ J $ attains its minimum at some $ \bar{u} \in u_0 + W_0^{1, 1}(\Omega) $.
	
	Furthermore, if $ (u,\xi) \mapsto L(x,u,\xi) $ is strictly convex for almost every $ x \in \Omega $, the minimizer is unique.
\end{theorem}

Applying Theorem \ref{T1}, we can also get the following proposition:

\begin{proposition}\label{Pr1}
	Let the functional $ J $ be defined on
	\[
		E = \{u \in W^{1, 1}(0,1) \: \vert \: u(0) = 0, \: u(1) = 1\}
	\]
	with the form
	\[
		J(u) = \int_0^1 \bigl[ \bigl(k \lvert \dot{u}(t) \rvert ^{p} + \ell \lvert u(t) \rvert^{p}\bigr)^{\frac{1}{p}} + c(t) \bigr] \mathrm{d}t
	\]
	for some $ k \in \mathbb{R}_{> 0} $, $ \ell \in \mathbb{R}_{\geq 0} $, $ c \in L^1(0,1) $ and $ p \in \mathbb{R}_{\geq 1} $. Then the infimum of $ J $ can be achieved at some $ \bar{u} \in E $ if and only if $ \ell = 0 $. 
	
	Moreover, if the infimum is attainable, the derivative $ \dot{\bar{u}}(t) $ is non-negative for almost every $ t \in (0,1) $.
\end{proposition}

Applying Proposition \ref{Pr1} we can directly answer the question in Example \ref{example1}, that is, $ J_1 $ cannot reach its infimum, but $ J_2 $ can. In fact, $ J_2 $ has infinitely many minimum points.

\bmhead{Outline}
The rest of the paper is organized as follows. In Section \ref{Preliminaries}, we introduce two important definitions and prove a technical lemma (Lemma \ref{Pweak}) that will be instrumental in the proof of Theorem \ref{T1}. In Section \ref{PMR}, we prove Theorem \ref{T1}, Theorem \ref{T2}, and Proposition \ref{Pr1}. In Appendix \ref{Poc}, we prove Corollary \ref{cor1} and Corollary \ref{cor2}. In Appendix \ref{Details}, we introduce more details about the syntactic category. In Appendix \ref{notations}, we provide a collection of notations and abbreviations used in this paper.

\section{Definitions and Lemmas}\label{Preliminaries}
To prove Theorem \ref{T1}, we introduce an operator $ \iota $ (see Definition \ref{definitionIs}) that embeds the Sobolev spaces $ W^{k,1} $ into some $ L^1 $ spaces. When considering a Sobolev space $ W^{k,1} $ defined on a non-empty open subset $ \Omega $ in a real Euclidean space $ \mathbb{R}^d $, we assign a replica $ \Omega_{\alpha} $ that is ``exactly the same" as the open set $ \Omega $ to each multi-index $ \alpha $, and take the (disjoint) union of these replicas $ \bigcup_{\lvert \alpha \rvert \leq k} \Omega_{\alpha} =: \Omega^{(k)} $ (see Definition \ref{definitionOmegak}) as the domain of the functions in the space $ L^1 $. We prove that the operator $ \iota $ is an isometry (see Lemma \ref{lemmaIs}). By converting it into a problem about the space $ L^1(\Omega^{(k)}) $, we can apply Theorem \ref{Dunford-Pettis} to complete the proof of Theorem \ref{T1}.

\begin{definition}[Construction of the disjoint union $ \Omega^{(k)} $]\textup{\cite[3.5 on page 61]{Adams}}\label{definitionOmegak}
	Let $ \Omega $ be a non-empty open set in $ \mathbb{R}^d $. For each $ \alpha $, let $ \Omega_{\alpha} $ be a different copy of $ \Omega $ lying in $ \mathbb{R}_{\alpha}^d $ which is a different copy of $ \mathbb{R}^d $ with respect to the multi-index $ \alpha $. Thus these $ \sum_{\lvert \alpha \rvert \leq k}  1 $ non-empty open sets are mutually disjoint. More mathematically, an isomorphism between $ \Omega $ and $ \Omega_{\alpha} $ refers to an isometry $ I_{\alpha} \colon \Omega \to \Omega_{\alpha} $ that also preserves the measure:
	\[
		\lambda_{\alpha}(\{I_{\alpha}(x) \: \vert \: x \in \omega\}) =: \lambda_{\alpha} \circ I_{\alpha}(\omega) = \lambda(\omega), \: \text{for all } \omega \subset \Omega\text{ with } \omega \text{ measurable}.
	\]
	Here, $\lambda_{\alpha}$ denotes the Lebesgue measure in $\mathbb{R}^d_\alpha$.
 
 The union of these $ \sum_{\lvert \alpha \rvert \leq k}  1 $ sets is denoted as
	\[
		\Omega^{(k)} := \bigsqcup_{\lvert \alpha \rvert \leq k} \Omega = \bigcup_{\lvert \alpha \rvert \leq k} \Omega_{\alpha}.
	\]
	The space $ \Omega^{(k)} $ naturally inherits the following structures from the Euclidean spaces:
	\begin{description}
		\item[\textbf{\textup{Measure:}}] A subset $\omega \subseteq \Omega^{(k)}$ is measurable if $\omega \cap \Omega_{\alpha}$ is Lebesgue measurable in $\mathbb{R}^d_\alpha$ for all $\alpha$. The measure of such $ \omega $ in $ \Omega^{(k)} $ is defined as the sum of the Lebesgue measures of $ \omega \cap \Omega_{\alpha} $ taken over $ \alpha $, i.e.
		\[
			\mu(\omega) := \sum_{\lvert \alpha \rvert \leq k}  \lambda_\alpha(\omega \cap \Omega_{\alpha}).
		\]
		\item[\textbf{\textup{Topology:}}] 
		The space $ \Omega^{(k)} $ is equipped with the disjoint union topology, i.e., a set $ \omega $ is open in $ \Omega^{(k)} $ if and only if $ \omega \cap \Omega_{\alpha} $ is open in $ \mathbb{R}_{\alpha}^d $ for every $ \alpha $.
	\end{description}
\end{definition}

Now we can define the isometry operator $ \iota $ from the Sobolev space $ W^{k,1}(\Omega) $ to the Lebesgue space $ L^1(\Omega^{(k)}) $. In this way, we transform the discussion on the compactness of a subset of $ W^{k,1}(\Omega) $ into the discussion on the compactness of its image in space $ L^1(\Omega^{(k)}) $.

\begin{definition}[Construction of the isometric operator $ \iota $]\textup{\cite[3.5 on page 61]{Adams}}\label{definitionIs}
	We construct an operator $ \iota $ from $ W^{k,1}( \Omega) $ to $ L^1(\Omega^{(k)}) $ as follows:
	\[
		\begin{split}
			\iota \colon \quad W^{k,1}(\Omega) &\to L^1(\Omega^{(k)})\\
			u &\mapsto \iota u,
		\end{split}
	\]
	where $\iota u: \Omega^{(k)} \to \mathbb{R}$ is given by
 \[\iota u(y) = (D^\alpha u) \circ I_{\alpha}^{-1}(y) \text{, if $y \in \Omega_\alpha$}.\]
	% \[
	% 	\begin{split}
	% 		\iota u \colon \quad \Omega^{(k)} &\to \mathbb{R}\\
	% 		y &\mapsto (\bigoplus_{\lvert \alpha \rvert \leq k} (D^{\alpha}u) \circ I_{\alpha}^{- 1})(y)\\
	% 		i.e. \: y &\mapsto (D^{\alpha}u) \circ I_{\alpha}^{- 1}(y), \: \text{if } y \in \Omega_{\alpha}.
	% 	\end{split}
	% \]
	That is, the restriction of function $ \iota u $ on each component of $ \Omega_{\alpha} $ in $\Omega^{(k)}$ is given by the function $ D^{\alpha}u $. 
 \end{definition}

 \begin{remark}
     From the perspective of category theory, let the objects in an indicator category $ I $ be the multi-indies, and the morphisms are only the identities. Then $ \Omega^{(k)} $ is the colimit of a $ G $ which maps each $ \alpha $ to $ \Omega_\alpha $. And $ L^1(\Omega^{(k)}) $ is the limit of a (contravariant) functor $ F $ which maps each $ \alpha $ to $ L^1(\Omega_\alpha) $. Therefore the operator $ \iota $ is the unique morphism from $ W^{k,1}( \Omega) $ to $ L^1(\Omega^{(k)}) $. See Chapter 1-3 of \cite{riehl2017category} for a precise explanation of these terminologies, and see the two illustrations in Figure~\ref{fig:colimit} and Figure~\ref{fig:L1}.
     
\begin{figure}[htb]
  \centering
  \begin{adjustbox}{width=0.5\textwidth}
  \begin{tikzpicture}
 \node[line width=0.5mm, minimum size=1cm] (A) at (0,5) {\Large $\alpha_1$};
\node[line width=0.5mm, minimum size=1cm] (B) at (0,4) {\Large $\alpha_2$};
\node[line width=0.5mm, minimum size=1cm] (C) at (0,3) {\Large $\vdots$};
\node[line width=0.5mm, minimum size=1cm] (D) at (0,2) {\Large $\alpha_n$};
\node[line width=0.5mm, minimum size=1cm] (E) at (0,1) {\Large $\vdots$};
\node[line width=0.5mm, minimum size=1cm] (F) at (0,0) {\normalsize $(\vert\alpha\vert \leq k)$};

 \node[line width=0.5mm, minimum size=1cm] (X) at (-1,3) {\Huge $ I:$};

  \node[line width=0.5mm, minimum size=1cm] (A1) at (3,5) {\Large $\Omega_{\alpha_1}$};
\node[line width=0.5mm, minimum size=1cm] (B1) at (3,4) {\Large $\Omega_{\alpha_2}$};
\node[line width=0.5mm, minimum size=1cm] (C1) at (3,3) {\Large $\vdots$};
\node[line width=0.5mm, minimum size=1cm] (D1) at (3,2) {\Large $\Omega_{\alpha_n}$};
\node[line width=0.5mm, minimum size=1cm] (E1) at (3,1) {\Large $\vdots$};

\draw[-{Latex[width=3mm]},->,decorate, decoration={zigzag}] (0.5, 2.9) to node [above=3pt] {$ \Omega_{\_} $} (2.5, 2.9);

\node[line width=0.5mm, minimum size=1cm] (A2) at (6.5,3) {\Large $\underrightarrow{\operatorname{colim}}\ \Omega_{\alpha_i} = \Omega^{(k)}$};

\draw[right hook-stealth] (A1) -- (A2);
\draw[right hook-stealth] (B1) -- (A2);
\draw[right hook-stealth] (C1) -- (A2);
\draw[right hook-stealth] (D1) -- (A2);
\draw[right hook-stealth] (E1) -- (A2);
\end{tikzpicture}
  \end{adjustbox}
  \caption{\centering Illustration of $ \Omega^{(k)} $ - the colimit.}
  \label{fig:colimit}
\end{figure}

\begin{figure}[htb]
  \centering
  \begin{adjustbox}{width=0.5\textwidth}
  \begin{tikzpicture}
 \node[line width=0.5mm, minimum size=1cm] (A) at (0,5) {\Large $\alpha_1$};
\node[line width=0.5mm, minimum size=1cm] (B) at (0,4) {\Large $\alpha_2$};
\node[line width=0.5mm, minimum size=1cm] (C) at (0,3) {\Large $\vdots$};
\node[line width=0.5mm, minimum size=1cm] (D) at (0,2) {\Large $\alpha_n$};
\node[line width=0.5mm, minimum size=1cm] (E) at (0,1) {\Large $\vdots$};
\node[line width=0.5mm, minimum size=1cm] (F) at (0,0) {\normalsize $(\vert\alpha\vert \leq k)$};

 \node[line width=0.5mm, minimum size=1cm] (X) at (-1,3) {\Huge $ I: $};

  \node[line width=0.5mm, minimum size=1cm] (A1) at (3,5) {\Large $L^1(\Omega_{\alpha_1})$};
\node[line width=0.5mm, minimum size=1cm] (B1) at (3,4) {\Large $L^1(\Omega_{\alpha_2})$};
\node[line width=0.5mm, minimum size=1cm] (C1) at (3,3) {\Large $\vdots$};
\node[line width=0.5mm, minimum size=1cm] (D1) at (3,2) {\Large $L^1(\Omega_{\alpha_n})$};
\node[line width=0.5mm, minimum size=1cm] (E1) at (3,1) {\Large $\vdots$};

\draw[-{Latex[width=3mm]},->,decorate, decoration={zigzag}] (0.5, 2.9) to node [above=3pt] {$ L^1(\Omega_{\_}) $} (2.5, 2.9);

\node[line width=0.5mm, minimum size=1cm] (A2) at (6,5) {\Large $W^{k,1}(\Omega)$};
\node[line width=0.5mm, minimum size=1cm] (B2) at (6,3) {\Large $\underleftarrow{\operatorname{lim}}\ L^1(\Omega_{\alpha_i})$};
\node[line width=0.5mm, minimum size=1cm] (B2) at (6,2.5) { $\vert \vert$};
\node[line width=0.5mm, minimum size=1cm] (B2) at (6,2) {\Large $L^1(\Omega^{(k)})$};

\draw [->] (A2) to node [above=3pt] {$D^{\alpha}$}  (A1);
\draw [->] (A2) -- (B1);
\draw [->] (A2) -- (D1);

\draw [->, dashed] (6, 4.5) to node [right=3pt] {$\exists ! \iota$} (6, 3.5);
\end{tikzpicture}
  \end{adjustbox}
  \caption{\centering Illustration of $ L^1(\Omega^{(k)}) $ - the limit.}
  \label{fig:L1}
\end{figure}
 \end{remark}

\begin{lemma}\label{lemmaIs}
The operator $\iota $ is a well-defined isometry. Let $ W $ be the range of $ \iota $. $W$ is a closed linear subspace of $L^1(\Omega^{(k)})$. Furthermore, $ W^{k, 1}(\Omega) $ and $ W $ are homeomorphic both in their norm and weak topologies.
\end{lemma}

\begin{proof}[Proof of \textup{Lemma~\ref{lemmaIs}}]
    We first check that $\iota $ is well-defined. Notice that for any $ y \in \Omega^{(k)} $, there must be only one multi-index $ \alpha $ such that $ y \in \Omega_{\alpha} $. Furthermore, for each $ \alpha $, the weak generalized partial derivative $ D^{\alpha}u $ is unique up to sets of measure zero in $ L^1(\Omega) $. Thus $ \iota $ is well-defined.

    Now we will check that $\iota $ is an isometry. Indeed, we have the following equalities
	\[
		\lVert \iota u \rVert_{L^1(\Omega^{(k)})} = \sum_{\lvert \alpha \rvert \leq k}  \lVert \iota u\vert_{\Omega_{\alpha}} \rVert_{L^1(\Omega_{\alpha})} = \sum_{\lvert \alpha \rvert \leq k}  \int_{\Omega}\lvert D^{\alpha} u \rvert \mathrm{d}\lambda=\lVert u \rVert_{{W^{k,1}(\Omega)}}.
	\]
 Since $\iota $ is an isometry between two complete metric spaces, the image of $\iota $ is also complete and thus closed in $L^1(\Omega^{(k)})$. Finally, it follows from \cite[see][Theorem 3.10 on page 61]{Brezis} that $ W^{k, 1}(\Omega) $ and $ W $ are homeomorphic both in their norm topologies and their weak topologies.
\end{proof}

Following this lemma, we denote the inverse operator of $\iota $ as $ \iota ^{- 1} \colon W \to W^{k, 1}(\Omega) $. As a useful result, we have the following lemma.

\begin{lemma}\label{Pweak}
	Suppose that (the same as Theorem \textup{\ref{T1}})
	\begin{description}
		\item[$ \mathbf{H^{k, 1}_{n.o.}}{[\Omega, \mathcal{F}]} $ \textup{(hypothesis):}] 
		$ \Omega $ is a non-empty open set in $ \mathbb{R}^d $, and $ \mathcal{F} $ is a subset of $ W^{k,1}(\Omega) $.
	\end{description}
	
	Then the following statements are equivalent:
	\begin{description}
		\item[$ \mathbf{C_w^{k, 1}}{[\Omega, \mathcal{F}]} $ \textup{(weak compactness in $ W^{k, 1} $):}] 
		$ \mathcal{F} $ is relatively weakly compact in $ W^{k,1}(\Omega) $.
		\item[$ \mathbf{C_w^1}{[\Omega^{(k)}, \iota \mathcal{F}]} $\textup{:}] $ \iota \mathcal{F} := \{\iota f \: \vert \: f \in \mathcal{F}\} $ is relatively weakly compact in $ L^1(\Omega^{(k)}) $.
	\end{description}
\end{lemma}

\begin{proof}[\textbf{Proof of \textup{Lemma \ref{Pweak}.}}]

We first need to prove that under the hypothesis $ \mathbf{H^{k, 1}_{n.o.}}{[\Omega, \mathcal{F}]} $, the following two topological spaces are homeomorphic. 
\begin{description}
    \item[a.] $ \bigl(W^{k,1}(\Omega), \sigma(W^{k,1}(\Omega), W^{k,1}(\Omega)^{\star})\bigr) $: the Sobolev space $ W^{k,1}(\Omega) $ equipped with the weak topology.
    \item[b.] $ \Bigl(W, \sigma\bigl(W, L^1(\Omega^{(k)})^{\star}\vert_W\bigr)\Bigr) $: the space $ W $ equipped with the subspace topology induced from the weak topology of the space $ L^1(\Omega^{(k)}) $, where $ L^1(\Omega^{(k)})^{\star}\vert_W := \{\phi \vert_W \: \vert \: \phi \in L^1(\Omega^{(k)})^{\star}\} $.
\end{description}

	Using Lemma \ref{lemmaIs}, we have
	\[
		\bigl(W^{k,1}(\Omega), \sigma(W^{k,1}(\Omega), W^{k,1}(\Omega)^{\star})\bigr) \cong \bigl(W, \sigma(W, W^{\star})\bigr).
	\]
	Next, we only need to show that
	\[
		\bigl(W, \sigma(W, W^{\star})\bigr) \cong \Bigl(W, \sigma\bigl(W, L^1(\Omega^{(k)})^{\star}\vert_W\bigr)\Bigr).
	\]
	
	Recall that the basis elements for the topology $ \sigma(W,W^{\star}) $ can be written as
	\begin{equation}\label{i1}
		\bigcap_{i \in F} \{w \in W \vert \: \lvert \left \langle \psi_i,w - w_i \right \rangle \rvert < r_i\},
	\end{equation}
	for some finite $ F $, and $ w_i \in W $, $ \psi_i \in W^{\star} $, $ r_i > 0 $.
	
	Using the analytic form of Hahn-Banach theorem \cite[see][Theorem 1.1 on page 1]{Brezis}, we have that for every $ \psi \in W^{\star} $, there exists $ \phi \in L^1(\Omega^{(k)})^{\star} $ such that
	\begin{align}
		&\left \langle \phi,w - w_i \right \rangle = \left \langle \psi,w - w_i \right \rangle, &&\text{for all } w \in W,\notag\\
		&\left \langle \phi,u - w_i \right \rangle \leq \lVert \psi \rVert_{W^{\star}}\lVert u \rVert_{L^1(\Omega^{(k)})}, &&\text{for all } u \in L^1(\Omega^{(k)}).\notag
	\end{align}
	Thus for each $ \psi_i $ in \eqref{i1}, there exists $ \phi_i \in L^1(\Omega^{(k)})^{\star} $ such that
	\[
		\{w \in W \vert \: \lvert \left \langle \psi_i,w - w_i \right \rangle \rvert < r_i\} = \{w \in W \vert \: \lvert \left \langle \phi_i,w - w_i \right \rangle \rvert < r_i\}
	\]
	It means that $ \sigma(W,W^{\star}) $ is stronger than $ \sigma(W, L^1(\Omega^{(k)})^{\star}\vert_W) $.
	
	Besides, similarly, every open set for $ \sigma(W, L^1(\Omega^{(k)})^{\star}) $ is a union of such sets:
	\begin{equation}\label{i2}
		\bigcap_{i \in F} \{w \in W \vert \: \lvert \left \langle \phi_i,w - w_i \right \rangle \rvert < r_i\},
	\end{equation}
	with $ F $ being a finite set, $ w_i \in W $, $ \phi_i \in L^1(\Omega^{(k)})^{\star} $, $ r_i > 0 $.
	
	Let $ \psi_i = \phi_i \vert_W \in W^{\star} $ for each $ \phi_i $ in \eqref{i2}. Then,
	\[
		\{w \in W \vert \: \lvert \left \langle \phi_i,w - w_i \right \rangle \rvert < r_i\} = \{w \in W \vert \: \lvert \left \langle \psi_i,w - w_i \right \rangle \rvert < r_i\}
	\]
	implies that $ \sigma(W, W^{\star}) $ is weaker than $ \sigma(W, L^1(\Omega^{(k)})^{\star}\vert_W) $.

	In summary, we obtain that $ \sigma(W, W^{\star}) $ and $ \sigma(W, L^1(\Omega^{(k)})^{\star}\vert_W) $ have the same topological basis. Thus
	\[
		\Bigl(W, \sigma\bigl(W, L^1(\Omega^{(k)})^{\star}\vert_W\bigr)\Bigr) \cong \bigl(W, \sigma(W, W^{\star})\bigr) \cong \bigl(W^{k,1}(\Omega), \sigma(W^{k,1}(\Omega), W^{k,1}(\Omega)^{\star})\bigr).
	\]
	
	We know from Lemma \ref{lemmaIs} that $ W $ is a closed linear subspace. Moreover, $ W $ is weakly closed by its convexity (linearity) using Mazur's lemma \cite[see][Theorem 2 on page 120]{Yosida}. It follows that relatively compact sets for $ \sigma(W, L^1(\Omega^{(k)})^{\star}\vert_W) $ coincide with relatively weakly compact sets in $ L^1(\Omega^{(k)}) $. Hence $ \mathcal{F} $ is relatively weakly compact in $ W^{k,1}(\Omega) $ if and only if $ \iota \mathcal{F} $ is relatively weakly compact in $ L^1(\Omega^{(k)}) $.
\end{proof}

\section{Proofs of Main Results}\label{PMR}

In this section we prove the main results in Section \ref{sec1}. We leave the proof of Corollary~\ref{cor1},~\ref{cor1'},~\ref{cor2}, and~\ref{cor2'} to Appendix~\ref{Poc}.

\begin{proof}[\textbf{Proof of \textup{Theorem \ref{T1}}.}]
	The idea is to show the following list of equivalences:
	\begin{align}
		&\mathbf{H^{k, 1}_{n.o.}} \wedge \mathbf{B^{k, 1}} \wedge \mathbf{EI^{k, 1}} \wedge \mathbf{E_{\infty}^{k, 1}}[\Omega, \mathcal{F}, \exists M, \exists \delta, \exists \omega]\label{i10}\\
		\Leftrightarrow &\mathbf{H^{k, 1}_{n.o.}}[\Omega, \mathcal{F}] \wedge \bigwedge_{\lvert \alpha \rvert \leq k} (\mathbf{H_{\boldsymbol{\sigma}}^1} \wedge \mathbf{B^1} \wedge \mathbf{EI^1} \wedge \mathbf{E_{\infty}^1}[\Omega_{\alpha}, \iota \mathcal{F}\vert_{\Omega_{\alpha}}, \exists M_{\alpha}, \exists \delta_{\alpha}, \exists \omega_{\alpha}])\label{i11}\\
		\Leftrightarrow &\mathbf{H^{k, 1}_{n.o.}}[\Omega, \mathcal{F}] \wedge (\mathbf{H_{\boldsymbol{\sigma}}^1} \wedge \mathbf{B^1} \wedge \mathbf{EI^1} \wedge \mathbf{E_{\infty}^1}[\Omega^{(k)}, \iota \mathcal{F}, \exists M, \exists \delta, \exists \omega^{(k)}])\label{i12}\\
		\Leftrightarrow &\mathbf{H^{k, 1}_{n.o.}}[\Omega, \mathcal{F}] \wedge (\mathbf{H_{\boldsymbol{\sigma}}^1} \wedge \mathbf{C_w^1}[\Omega^{(k)}, \iota \mathcal{F}])\label{i13}\\
		\Leftrightarrow &\mathbf{H^{k, 1}_{n.o.}} \wedge \mathbf{C_w^{k, 1}}[\Omega, \mathcal{F}]\label{i14}
	\end{align}
    For convenience of the reader, we summarize the main idea behind proving the equivalences of these five statements in Figure~\ref{fig:proof}.\\

    \begin{figure}[htb]
  \centering
  \begin{adjustbox}{width=\textwidth}
  \begin{tikzpicture}

% -----------------------------------------
 \node[line width=0.5mm, minimum size=1cm] (X1) at (0,5) {\Huge $\bigwedge$:};

 \node[line width=0.5mm, minimum size=1cm] (X2) at (3,5.3) {\Large $\mathbf{H^{k,1}_{n.o.}}$};
 
\node[line width=0.5mm, minimum size=1cm] (B1) at (0,4) {\Large $\mathbf{H^1_\sigma[\bullet_{(0, ..., 0)}]}$};
\node[line width=0.5mm, minimum size=1cm] (C) at (0,3) {\Large $\mathbf{B^1[\bullet_{(0, ..., 0)}]}$};
\node[line width=0.5mm, minimum size=1cm] (D) at (0,2) {\Large $\mathbf{EI^1[\bullet_{(0,...,0)}]}$};
\node[line width=0.5mm, minimum size=1cm] (E) at (0,1) {\Large $\mathbf{E^1_\infty[\bullet_{(0,...,0)}]}$};

\node[line width=0.5mm, minimum size=1cm] (B) at (1.7,3.8) {\Large $...$};
\node[line width=0.5mm, minimum size=1cm] (B) at (1.7,2.8) {\Large $...$};
\node[line width=0.5mm, minimum size=1cm] (B) at (1.7,1.8) {\Large $...$};
\node[line width=0.5mm, minimum size=1cm] (B) at (1.7,0.8) {\Large $...$};

\node[line width=0.5mm, minimum size=1cm] (B2) at (3,4) {\Large $\mathbf{H^1_\sigma[\bullet_{\alpha}]}$};
\node[line width=0.5mm, minimum size=1cm] (C) at (3,3) {\Large $\mathbf{B^1[\bullet_{\alpha}]}$};
\node[line width=0.5mm, minimum size=1cm] (D) at (3,2) {\Large $\mathbf{EI^1[\bullet_{\alpha}]}$};
\node[line width=0.5mm, minimum size=1cm] (E) at (3,1) {\Large $\mathbf{E^1_\infty[\bullet_{\alpha}]}$};

\node[line width=0.5mm, minimum size=1cm] (B) at (4.25,3.8) {\Large $...$};
\node[line width=0.5mm, minimum size=1cm] (B) at (4.25,2.8) {\Large $...$};
\node[line width=0.5mm, minimum size=1cm] (B) at (4.25,1.8) {\Large $...$};
\node[line width=0.5mm, minimum size=1cm] (B) at (4.25,0.8) {\Large $...$};

\node[line width=0.5mm, minimum size=1cm] (B3) at (6,4) {\Large $\mathbf{H^1_\sigma[\bullet_{(0, ..., k)}]}$};
\node[line width=0.5mm, minimum size=1cm] (C) at (6,3) {\Large $\mathbf{B^1[\bullet_{(0, ..., k)}]}$};
\node[line width=0.5mm, minimum size=1cm] (D) at (6,2) {\Large $\mathbf{EI^1[\bullet_{(0, ..., k)}]}$};
\node[line width=0.5mm, minimum size=1cm] (E) at (6,1) {\Large $\mathbf{E^1_\infty[\bullet_{(0, ...,k)}]}$};

\draw[double,->] (X2) -- (B1);
\draw[double,->] (X2) -- (B2);
\draw[double,->] (X2) -- (B3);
% -----------------------------------------

\node[line width=0.5mm, minimum size=1cm] (B3) at (1,-0.5) {\Large $\bigwedge:$};
\node[line width=0.5mm, minimum size=1cm] (B3) at (1,-1.8) {\Large $\mathbf{H^{k,1}_{n.o.}}$};
\node[line width=0.5mm, minimum size=1cm] (B3) at (1,-2.8) {\Large $\mathbf{C^{k,1}_W}$};
% -----------------------------------------

\node[line width=0.5mm, minimum size=1cm] (B4) at (5,-0.5) {\Large $\bigwedge:\mathbf{H^{k,1}_{n.o.}}$};
\node[line width=0.5mm, minimum size=1cm] (B5) at (5,-1.8) {\Large $\mathbf{H^1_\sigma[\bullet^{(k)}]}$};
\node[line width=0.5mm, minimum size=1cm] (B6) at (5,-2.8) {\Large $\mathbf{C^1_W[\bullet^{(k)}]}$};
\draw[double,->] (B4) -- (B5);

% -----------------------------------------

\node[line width=0.5mm, minimum size=1cm] (B4) at (10,3) {\Large $\bigwedge:\mathbf{H^{k,1}_{n.o.}}$};
\node[line width=0.5mm, minimum size=1cm] (B5) at (10,1.7) {\Large $\mathbf{H^1_\sigma[\bullet^{(k)}]}$};
\node[line width=0.5mm, minimum size=1cm] (B6) at (10,0.7) {\Large $\mathbf{B^1[\bullet^{(k)}]}$};
\node[line width=0.5mm, minimum size=1cm] (B7) at (10,-0.3) {\Large $\mathbf{EI^1[\bullet^{(k)}]}$};
\node[line width=0.5mm, minimum size=1cm] (B8) at (10,-1.3) {\Large $\mathbf{E^1_\infty[\bullet^{(k)}]}$};
\draw[double,->] (B4) -- (B5);
% -----------------------------------------
\node[line width=0.5mm, minimum size=1cm] (B4) at (-4,3) {\Large $\bigwedge:$};
\node[line width=0.5mm, minimum size=1cm] (B5) at (-4,1.7) {\Large $\mathbf{H^{k,1}_{n.o.}}$};
\node[line width=0.5mm, minimum size=1cm] (B6) at (-4,0.7) {\Large $\mathbf{B^{k,1}}$};
\node[line width=0.5mm, minimum size=1cm] (B7) at (-4,-0.3) {\Large $\mathbf{EI^{k,1}}$};
\node[line width=0.5mm, minimum size=1cm] (B8) at (-4,-1.3) {\Large $\mathbf{E^{k,1}_\infty}$};

% --------------------------------------------
% Drawing Rectangles
\draw[rounded corners=4.4mm, line width=0.5mm] (-1.5, 5.7) rectangle (7.5, 0.5);
\draw[rounded corners=4.4mm, line width=0.5mm] (0, 0) rectangle (2, -3.5);
\draw[rounded corners=4.4mm, line width=0.5mm] (3.8, 0) rectangle (6.2, -3.5);
\draw[rounded corners=4.4mm, line width=0.5mm] (8.8, 3.5) rectangle (11.2, -2);
\draw[rounded corners=4.4mm, line width=0.5mm] (-5, 3.5) rectangle (-3, -2);

% --------------------------------------------
\draw[double,<->] (-1.6,3) to node [midway,above,sloped] {\tiny Def. of $\iota$} (-2.8,0.5);

\draw[double,<->] (-2.8,0) to node [above,sloped, align=center] {\tiny Weak Compactness\\\tiny Criterion in $W^{k,1}$} (0,-1.6);

\draw[double,<->] (2, -1.6) to node [above,sloped, align=center] {\tiny Lemma \ref{Pweak}} (3.8,-1.6);

\draw[double,<->] (8.8,0) to node [above,sloped, align=center] {\tiny Dunford-Pettis Theorem} (6.2,-1.6);

\draw[double,<->] (7.6,3) to node [midway,above,sloped] {\tiny Def. of $\Omega^{(k)}$} (8.8,0.5);

\end{tikzpicture}
  \end{adjustbox}
  \caption{\centering Summary of the Proof of Theorem \ref{T1}.}
  \label{fig:proof}
\end{figure}

	Now we move on to the actual proof. As a result that will be used repeatedly, we first prove that:
	\begin{equation}\label{i01}
		\begin{split}
			\mathbf{H^{k, 1}_{n.o.}}[\Omega, \mathcal{F}] &\Leftrightarrow \mathbf{H^{k, 1}_{n.o.}}[\Omega, \mathcal{F}] \wedge \bigwedge_{\lvert \alpha \rvert \leq k} (\mathbf{H_{\boldsymbol{\sigma}}^1}[\Omega_{\alpha}, \iota \mathcal{F}\vert_{\Omega_{\alpha}}])\\
			&\Leftrightarrow \mathbf{H^{k, 1}_{n.o.}}[\Omega, \mathcal{F}] \wedge \mathbf{H_{\boldsymbol{\sigma}}^1}[\Omega^{(k)}, \iota \mathcal{F}].
		\end{split}
	\end{equation}
	In fact, since $ \Omega_{\alpha} $ is a $ \sigma $-finite measure space (as an open set in $ \mathbb{R}_{\alpha}^d $ which is isometric with $ \mathbb{R}^d $), we have $ \mathbf{H^{k, 1}_{n.o.}}[\Omega, \mathcal{F}] \Rightarrow \mathbf{H_{\boldsymbol{\sigma}}^1}[\Omega_{\alpha}, \iota \mathcal{F}\vert_{\Omega_{\alpha}}] $ for every $ \alpha $. By the definition of the general product in $ \mathfrak{C} $ (see Definition \ref{category}), we immediately obtain that:
	\[
		\mathbf{H^{k, 1}_{n.o.}}[\Omega, \mathcal{F}] \Rightarrow \mathbf{H^{k, 1}_{n.o.}}[\Omega, \mathcal{F}] \wedge \bigwedge_{\lvert \alpha \rvert \leq k} (\mathbf{H_{\boldsymbol{\sigma}}^1}[\Omega_{\alpha}, \iota \mathcal{F}\vert_{\Omega_{\alpha}}]) \Rightarrow \mathbf{H^{k, 1}_{n.o.}}[\Omega, \mathcal{F}].
	\]
	Moreover, since $ \Omega^{(k)} $ is also a $ \sigma $-finite measure space, we have
	\[
		\bigwedge_{\lvert \alpha \rvert \leq k} (\mathbf{H_{\boldsymbol{\sigma}}^1}[\Omega_{\alpha}, \iota \mathcal{F}\vert_{\Omega_{\alpha}}]) \Leftrightarrow \mathbf{H_{\boldsymbol{\sigma}}^1}[\Omega^{(k)}, \iota \mathcal{F}].
	\]
	Thus we get \eqref{i01}. Herein we show \eqref{i10} $ \Leftrightarrow $ \eqref{i11} $ \Leftrightarrow \cdots \Leftrightarrow $ \eqref{i14} one by one. These are all simple verifications.

	\noindent Proof of ``\eqref{i10} $ \Rightarrow $ \eqref{i11}'':
	
	Suppose $ \mathbf{H^{k, 1}_{n.o.}} \wedge \mathbf{B^{k, 1}} \wedge \mathbf{EI^{k, 1}} \wedge \mathbf{E_{\infty}^{k, 1}}[\Omega, \mathcal{F}, \exists M, \exists \delta, \exists \omega] $, and let
	\[
		M_{\alpha} = M, \: \delta_{\alpha} = \delta, \: \omega_{\alpha} = I_{\alpha} \circ \omega.
	\]
	It is not difficult to verify that:
	\begin{itemize}
		\item[1.] For any $ \iota f\vert_{\Omega_{\alpha}} \in \iota \mathcal{F}\vert_{\Omega_{\alpha}} := \{(\iota f)\vert_{\Omega_{\alpha}} \: \vert \: f \in \mathcal{F}\} $, its norm satisfies
		\[
			\lVert \iota f\vert_{\Omega_{\alpha}} \rVert_{L^1(\Omega_{\alpha})} = \lVert D^{\alpha}f \rVert_{L^1(\Omega)} \leq \sum_{\lvert \beta \rvert \leq k}  \lVert D^{\beta}f \rVert_{L^1(\Omega)} = \lVert f \rVert_{W^{k, 1}(\Omega)} \leq M = M_{\alpha},
		\]
		which implies $ \mathbf{B^1}[\Omega_{\alpha}, \iota \mathcal{F}\vert_{\Omega_{\alpha}}, \exists M_{\alpha}] $ for every $ \alpha $.
		\item[2.] Similarly, for each $ \iota f\vert_{\Omega_{\alpha}} \in \iota \mathcal{F}\vert_{\Omega_{\alpha}}$, for any $ \varepsilon \in \mathbb{R}_{> 0} $ and for every measurable set $ A \subseteq \Omega_{\alpha} $ with its measure $ \lambda_{\alpha}(A) < \delta_{\alpha}(\varepsilon) $, since $ \lambda \circ I_{\alpha}^{- 1}(A) = \lambda(A) < \delta_{\alpha}(\varepsilon) = \delta(\varepsilon) $, we have
		\[
			\int_A \lvert \iota f\vert_{\Omega_{\alpha}} \rvert \mathrm{d}\lambda_{\alpha} = \int_{I_{\alpha}^{- 1}(A)} \lvert D^{\alpha}f \rvert \mathrm{d}\lambda \leq \sum_{\lvert \beta \rvert \leq k}  \int_{I_{\beta}^{- 1}(A)} \lvert D^{\beta}f \rvert \mathrm{d}\lambda < \varepsilon.
		\]
		Thus we get $ \mathbf{EI^1}[\Omega_{\alpha}, \iota \mathcal{F}\vert_{\Omega_{\alpha}}, \exists \delta_{\alpha}] $ for every $ \alpha $.
		\item[3.] Again, for each $ \iota f\vert_{\Omega_{\alpha}} \in \iota \mathcal{F}\vert_{\Omega_{\alpha}} $ and for any $ \varepsilon \in \mathbb{R}_{> 0} $, the inequality
		\[
			\int_{\Omega_{\alpha} \setminus \omega_{\alpha}(\varepsilon)} \lvert \iota f\vert_{\Omega_{\alpha}} \rvert \mathrm{d}\lambda_{\alpha} = \int_{\Omega \setminus \omega(\varepsilon)} \lvert D^{\alpha}f \rvert \mathrm{d}\lambda \leq \sum_{\lvert \beta \rvert \leq k}  \int_{\Omega \setminus \omega(\varepsilon)} \lvert D^{\beta}f \rvert \mathrm{d}\lambda < \varepsilon
		\]
		implies $ \mathbf{E_{\infty}^1}[\Omega_{\alpha}, \iota \mathcal{F}\vert_{\Omega_{\alpha}}, \exists \omega_{\alpha}] $ for each $ \alpha $.
	\end{itemize}
	Hence by the definition of general products, there exists the morphism ``\eqref{i10} $ \Rightarrow $ \eqref{i11}''.

	\noindent Proof of ``\eqref{i11} $ \Rightarrow $ \eqref{i10}'':

	Suppose $ \mathbf{H^{k, 1}_{n.o.}}[\Omega, \mathcal{F}] \wedge \bigwedge_{\lvert \alpha \rvert \leq k} (\mathbf{H_{\boldsymbol{\sigma}}^1} \wedge \mathbf{B^1} \wedge \mathbf{EI^1} \wedge \mathbf{E_{\infty}^1}[\Omega_{\alpha}, \iota \mathcal{F}\vert_{\Omega_{\alpha}}, \exists M_{\alpha}, \exists \delta_{\alpha}, \exists \omega_{\alpha}]) $, and let
	\[
		M = \sum_{\lvert \alpha \rvert \leq k}  M_{\alpha}, \: \delta \colon \varepsilon \mapsto \min_{\lvert \alpha \rvert \leq k} \delta_{\alpha}(\frac{\varepsilon}{\sum_{\lvert \alpha \rvert \leq k}  1}), \: \omega \colon \varepsilon \mapsto \bigcup_{\lvert \alpha \rvert \leq k} I_{\alpha}^{- 1} \circ \omega_{\alpha}(\frac{\varepsilon}{\sum_{\lvert \alpha \rvert \leq k}  1}).
	\]
	We obtain that:
	\begin{itemize}
		\item[1.] For any $ f \in \mathcal{F} $, its norm satisfies that
		\[
			\lVert f \rVert_{W^{k, 1}} = \sum_{\lvert \alpha \rvert \leq k}  \lVert D^{\alpha}f \rVert_{L^1(\Omega)} = \sum_{\lvert \alpha \rvert \leq k}  \lVert \iota f\vert_{\Omega_{\alpha}} \rVert_{L^1(\Omega_{\alpha})} \leq \sum_{\lvert \alpha \rvert \leq k}  M_{\alpha} = M,
		\]
		which implies $ \mathbf{B^{k, 1}}[\Omega, \mathcal{F}, \exists M] $.
		\item[2.] For each $ f \in \mathcal{F} $, for any $ \varepsilon \in \mathbb{R}_{> 0} $ and for every measurable set $ A \subseteq \Omega $ with its measure $ \lambda(A) < \delta(\varepsilon) $, since $ \lambda_{\alpha} \circ I_{\alpha}(A) = \lambda(A) < \delta(\varepsilon) = \min_{\lvert \alpha \rvert \leq k} \delta_{\alpha}(\frac{\varepsilon}{\sum_{\lvert \alpha \rvert \leq k}  1}) \leq \delta_{\alpha}(\frac{\varepsilon}{\sum_{\lvert \alpha \rvert \leq k}  1}) $ for all $ \alpha $, we have
		\[
			\sum_{\lvert \alpha \rvert \leq k}  \int_A \lvert D^{\alpha}f \rvert \mathrm{d}\lambda = \sum_{\lvert \alpha \rvert \leq k}  \int_{I_{\alpha}(A)} \lvert \iota f\vert_{\Omega_{\alpha}} \rvert \mathrm{d}\lambda_{\alpha} < \sum_{\lvert \alpha \rvert \leq k}  \frac{\varepsilon}{\sum_{\lvert \alpha \rvert \leq k}  1} = \varepsilon.
		\]
		Thus we get $ \mathbf{EI^{k, 1}}[\Omega, \mathcal{F}, \exists \delta] $.
		\item[3.] In the same way, for each $ f \in \mathcal{F} $ and any $ \varepsilon \in \mathbb{R}_{> 0} $, the inequality
		\[
			\sum_{\lvert \alpha \rvert \leq k}  \int_{\Omega \setminus \omega(\varepsilon)} \lvert D^{\alpha}f \rvert \mathrm{d}\lambda \leq \sum_{\lvert \alpha \rvert \leq k}  \int_{\Omega_{\alpha} \setminus \omega_{\alpha}(\frac{\varepsilon}{\sum_{\lvert \alpha \rvert \leq k}  1})} \lvert \iota f\vert_{\Omega_{\alpha}} \rvert \mathrm{d}\lambda_{\alpha} < \sum_{\lvert \alpha \rvert \leq k}  \frac{\varepsilon}{\sum_{\lvert \alpha \rvert \leq k}  1} = \varepsilon
		\]
		implies $ \mathbf{E_{\infty}^{k, 1}}[\Omega, \mathcal{F}, \exists \omega] $.
	\end{itemize}

	\noindent Proof of ``\eqref{i11} $ \Rightarrow $ \eqref{i12}'':

	Suppose $ \mathbf{H^{k, 1}_{n.o.}}[\Omega, \mathcal{F}] \wedge \bigwedge_{\lvert \alpha \rvert \leq k} (\mathbf{H_{\boldsymbol{\sigma}}^1} \wedge \mathbf{B^1} \wedge \mathbf{EI^1} \wedge \mathbf{E_{\infty}^1}[\Omega_{\alpha}, \iota \mathcal{F}\vert_{\Omega_{\alpha}}, \exists M_{\alpha}, \exists \delta_{\alpha}, \exists \omega_{\alpha}]) $, and let
	\[
		M = \sum_{\lvert \alpha \rvert \leq k}  M_{\alpha}, \: \delta \colon \varepsilon \mapsto \min_{\lvert \alpha \rvert \leq k} \delta_{\alpha}(\frac{\varepsilon}{\sum_{\lvert \alpha \rvert \leq k}  1}), \: \omega^{(k)} \colon \varepsilon \mapsto \bigcup_{\lvert \alpha \rvert \leq k} \omega_{\alpha}(\frac{\varepsilon}{\sum_{\lvert \alpha \rvert \leq k}  1}).
	\]
	Almost the same as the previous argument, we can verify that:
	\begin{itemize}
		\item[1.] For every $ \iota f \in \iota \mathcal{F} $, its norm satisfies
		\[
			\lVert \iota f \rVert_{L^1(\Omega^{(k)})} = \sum_{\lvert \alpha \rvert \leq k}  \lVert \iota f\vert_{\Omega_{\alpha}} \rVert_{L^1(\Omega_{\alpha})} \leq \sum_{\lvert \alpha \rvert \leq k}  M_{\alpha} = M,
		\]
		which implies $ \mathbf{B^1}[\Omega^{(k)}, \iota \mathcal{F}, \exists M] $.
		\item[2.] For each $ \iota f \in \iota \mathcal{F} $, for any $ \varepsilon \in \mathbb{R}_{> 0} $, and for every measurable set $ A \subseteq \Omega^{(k)} $ with its measure $ \mu(A) < \delta(\varepsilon) $, since $ \lambda_{\alpha}(A \cap \Omega_{\alpha}) \leq \mu(A) < \delta(\varepsilon) = \min_{\lvert \alpha \rvert \leq k} \delta_{\alpha}(\frac{\varepsilon}{\sum_{\lvert \alpha \rvert \leq k}  1}) \leq \delta_{\alpha}(\frac{\varepsilon}{\sum_{\lvert \alpha \rvert \leq k}  1}) $ for all $ \alpha $, we have
		\[
			\int_A \lvert \iota f \rvert \mathrm{d}\mu = \sum_{\lvert \alpha \rvert \leq k}  \int_{A \cap \Omega_{\alpha}} \lvert \iota f\vert_{\Omega_{\alpha}} \rvert \mathrm{d}\lambda_{\alpha} < \sum_{\lvert \alpha \rvert \leq k}  \frac{\varepsilon}{\sum_{\lvert \alpha \rvert \leq k}  1} = \varepsilon.
		\]
		Thus we get $ \mathbf{EI^1}[\Omega^{(k)}, \iota \mathcal{F}, \exists \delta] $.
		\item[3.] For each $ \iota f \in \iota \mathcal{F} $ and any $ \varepsilon \in \mathbb{R}_{> 0} $, we have the inequality
		\[
			\int_{\Omega \setminus \omega(\varepsilon)} \lvert \iota f \rvert \mathrm{d}\mu = \sum_{\lvert \alpha \rvert \leq k}  \int_{\Omega_{\alpha} \setminus \omega_{\alpha}(\frac{\varepsilon}{\sum_{\lvert \alpha \rvert \leq k}  1})} \lvert \iota f\vert_{\Omega_{\alpha}} \rvert \mathrm{d}\lambda_{\alpha} < \sum_{\lvert \alpha \rvert \leq k}  \frac{\varepsilon}{\sum_{\lvert \alpha \rvert \leq k}  1} = \varepsilon,
		\]
		which leads to $ \mathbf{E_{\infty}^1}[\Omega^{(k)}, \iota \mathcal{F}, \exists \omega^{(k)}] $ immediately.
	\end{itemize}

	Once again, by the definition of general products, we have ``\eqref{i11} $ \Rightarrow $ \eqref{i12}''.

	\noindent Proof of ``\eqref{i12} $ \Rightarrow $ \eqref{i11}'':

	Suppose $ \mathbf{H^{k, 1}_{n.o.}}[\Omega, \mathcal{F}] \wedge (\mathbf{H_{\boldsymbol{\sigma}}^1} \wedge \mathbf{B^1} \wedge \mathbf{EI^1} \wedge \mathbf{E_{\infty}^1}[\Omega^{(k)}, \iota \mathcal{F}, \exists M, \exists \delta, \exists \omega^{(k)}]) $, and let
	\[
		M_{\alpha} = M, \: \delta_{\alpha} = \delta, \: \omega_{\alpha} \colon \varepsilon \mapsto \omega^{(k)}(\varepsilon) \cap \Omega_{\alpha}.
	\]
	We obtain that:
	\begin{itemize}
		\item[1.] For every $ \iota f\vert_{\Omega_{\alpha}} \in \iota \mathcal{F}\vert_{\Omega_{\alpha}} $, its norm satisfies
		\[
			\lVert \iota f\vert_{\Omega_{\alpha}} \rVert_{L^1(\Omega_{\alpha})} \leq \lVert \iota f \rVert_{L^1(\Omega^{(k)})} \leq M = M_{\alpha},
		\]
		which implies $ \mathbf{B^1}[\Omega_{\alpha}, \iota \mathcal{F}\vert_{\Omega_{\alpha}}, \exists M_{\alpha}] $ for every $ \alpha $.
		\item[2.] For each $ \iota f\vert_{\Omega_{\alpha}} \in \iota \mathcal{F}\vert_{\Omega_{\alpha}} $, for any $ \varepsilon \in \mathbb{R}_{> 0} $, and for every measurable set $ A \subseteq \Omega_{\alpha} $ with its measure $ \lambda_{\alpha}(A) < \delta_{\alpha}(\varepsilon) $, we have
		\[
			\int_A \lvert \iota f\vert_{\Omega_{\alpha}} \rvert \mathrm{d}\lambda_{\alpha} = \int_A \lvert \iota f \rvert \mathrm{d}\mu < \varepsilon.
		\]
		Thus we get $ \mathbf{EI^1}[\Omega_{\alpha}, \iota \mathcal{F}\vert_{\Omega_{\alpha}}, \exists \delta_{\alpha}] $ for every $ \alpha $.
		\item[3.] For each $ \iota f\vert_{\Omega_{\alpha}} \in \iota \mathcal{F}\vert_{\Omega_{\alpha}} $, the inequality
		\[
			\int_{\Omega_{\alpha} \setminus \omega_{\alpha}(\varepsilon)} \lvert \iota f\vert_{\Omega_{\alpha}} \rvert \mathrm{d}\lambda_{\alpha} \leq \int_{\Omega^{(k)} \setminus \omega^{(k)}(\varepsilon)} \lvert \iota f \rvert \mathrm{d}\mu < \varepsilon
		\]
		implies $ \mathbf{E_{\infty}^1}[\Omega_{\alpha}, \iota \mathcal{F}\vert_{\Omega_{\alpha}}, \exists \omega_{\alpha}] $ for every $ \alpha $.
	\end{itemize}
	
	These arguments complete the proof of ``\eqref{i12} $ \Rightarrow $ \eqref{i11}''.

	Isomorphism ``\eqref{i12} $ \Leftrightarrow $ \eqref{i13}'' can be obtained immediately from Theorem \ref{Dunford-Pettis}. Using Lemma \ref{Pweak} and \eqref{i01}, isomorphism ``\eqref{i13} $ \Leftrightarrow $ \eqref{i14}'' is obvious.

	Finally, the proof of Theorem \ref{T1} is completed by compositing.
\end{proof}

\begin{proof}[\textbf{Proof of \textup{Theorem \ref{T1'}.}}]
    We only need to prove the following two.
    \begin{description}
        \item[\textup{1.}] $ \mathbf{H^{k, 1}_{b.o.cone}} \wedge \mathbf{B^{\lvert \alpha \rvert \in \{0, k\},1}}{[\Omega, \mathcal{F}, \exists M]} \Rightarrow \mathbf{H^{k, 1}_{b.o.cone}} \wedge \mathbf{B^{k,1}}{[\Omega, \mathcal{F}, \exists M']} $. 
    \end{description}
    
    It can be obtained directly from the interpolation inequality \cite[see][Theorem 5.2 on page 135]{Adams}.
    
    \begin{description}  
        \item[\textup{2.}] $ \mathbf{H^{k, 1}_{b.o.cone}}\wedge \mathbf{B^{\lvert \alpha \rvert \in \{0, k\},1}} \wedge \mathbf{EI^{\lvert \alpha \rvert = k,1}}{[\Omega, \mathcal{F}, \exists M, \exists \delta]} \Rightarrow \mathbf{H^{k, 1}_{b.o.cone}}\wedge \mathbf{B^{\lvert \alpha \rvert \in \{0, k\},1}} \wedge \mathbf{EI^{k,1}}{[\Omega, \mathcal{F}, \exists M', \exists \delta']} $.
    \end{description}

        From the above discussion we already know that $ \mathcal{F} $ is bounded in $ W^{k,1}(\Omega) $. Applying the Rellich-Kondrachov theorem \cite[see][Theorem 6.3 on page 168]{Adams}, we obtain that the embedding from $ W^{k,1}(\Omega) $ to $ W^{k-1,1}(\Omega) $ is compact (When $ k = 1 $, the notation $ W^{0,1} $ is regarded as $ L^1 $). Hence $ \iota\mathcal{F} $ is precompact in $ L^1(\Omega^{(k)}) $ and obviously, it is relatively weakly compact. By Lemma \ref{Pweak}, we have that $ \mathcal{F} $ is relatively weakly compact in $ W^{k,1}(\Omega) $. Then applying Theorem \ref{T1}, the statement $ \mathbf{EI^{k,1}}{[\Omega, \mathcal{F}, \exists \delta']} $ is satisfied.
        
    Thus the proof is completed.
\end{proof}

\begin{proof}[\textbf{Proof of \textup{Theorem \ref{T2}.}}]
	Denote
	\[
		W_{u_0}^{1,1}(\Omega) = u_0 + W_0^{1, 1}(\Omega) \text{ and } m = \inf_{u \in W_{u_0}^{1,1}(\Omega)} J(u),
	\]
	and observe that $ m \leq J(u_0) < + \infty $. Integrating both sides of \eqref{E4.1} with respect to $ x $, we have
	\begin{equation}\label{E4.4}
		\begin{split}
			J(u) &= \int_{\Omega} L(x,u,\nabla u) \mathrm{d}x \geq a \int_{\Omega} \lvert \nabla  u \rvert \mathrm{d}x + \int_{\Omega} c(x) \mathrm{d}x\\
			&\geq a \lVert \nabla u \rVert_{L^1(\Omega)} - \lVert c \rVert_{L^1(\Omega)} > - \infty.
		\end{split}
	\end{equation}
	
	Next, we show that $ m > - \infty $. Let $ u = u_0 + v $ for some $ v \in W_0^{1, 1}(\Omega) $. Invoking Poincar\'{e}'s inequality \cite[see][Corollary 9.19 on page 290]{Brezis}, there exists a positive constant $ C $ such that
	\[
		\lVert \nabla v \rVert_{L^1(\Omega)} \geq C \lVert  v \rVert_{W_0^{1, 1}(\Omega)}.
	\]
	Combining with the triangle inequality, we have
	\[
		\begin{split}
			\lVert \nabla u \rVert_{L^1(\Omega)} + \lVert u_0 \rVert_{W^{1, 1}(\Omega)} &\geq \lVert \nabla u \rVert_{L^1(\Omega)} + \lVert \nabla u_0 \rVert_{L^1(\Omega)}\\
			&\geq \lVert \nabla u-\nabla u_0 \rVert_{L^1(\Omega)}\\
			&\geq C \lVert u - u_0 \rVert_{W_0^{1, 1}(\Omega)}\\
			&= C \lVert u- u_0 \rVert_{W^{1, 1}(\Omega)}\\
			&\geq C \lVert u \rVert_{W^{1, 1}(\Omega)} - C \lVert u_0 \rVert_{W^{1, 1}(\Omega)}.
		\end{split}
	\]
	In short,
	\begin{equation}\label{E4.7}
		\lVert \nabla u \rVert_{L^1(\Omega)} \geq C \lVert u \rVert_{W^{1, 1}(\Omega)} - (C + 1) \lVert u_0 \rVert_{W^{1, 1}(\Omega)}.
	\end{equation}
	Denote
	\[
		\gamma(a,C,u_0,c) = a(C+1) \lVert u_0 \rVert_{W^{1, 1}(\Omega)} + \lVert c \rVert_{L^1(\Omega)}.
	\]
	According to \eqref{E4.4} and \eqref{E4.7}, we obtain that
	\begin{equation}\label{E4.8}
		J(u) \geq aC \lVert u \rVert_{W^{1, 1}(\Omega)} - \gamma \geq - \gamma, \: \forall u \in W_{u_0}^{1,1}(\Omega).
	\end{equation}
	Thus we have proved that
	\[
		m = \inf_{u \in W_{u_0}^{1,1}(\Omega)} J(u) \geq - \gamma > - \infty.
	\]
	
	Hence, for an arbitrary minimizing sequence $ \{u_n\} \subset W_{u_0}^{1,1}(\Omega) $ and then for some integer $ n $ sufficiently large, we obtain that
	\[
		m + 1 \geq J(u_n).
	\]
	Combining with \eqref{E4.8}, we have
	\[
		m + 1 \geq aC \lVert u_n \rVert_{W^{1, 1}(\Omega)} - \gamma.
	\]
	It means that there exists a number $ M = M(m,a,C,\gamma) < + \infty $ such that
	\[
		\lVert u_n \rVert_{W^{1, 1}(\Omega)} \leq M.
	\]
	
	If we could find some suitable sequence $ \{u_n\} $ turning out to satisfies the statement $ \mathbf{EI^{\lvert \alpha \rvert = 1,1}} $ in Theorem \ref{T1'}, we should have that $ \{u_n\} $ is relatively weakly compact. Thus there exists a subsequence of $ \{u_n\} $ (also denoted as $ \{u_n\}$ for simplicity) converging weakly to some $ \bar{u} \in W_{u_0}^{1,1}(\Omega) $.
	
	Using Theorem 3.23 in \cite[page 96]{Dacorogna} to obtain that functional $ J $ is weakly lower semicontinuous, we have
	\[
		\liminf_{n \to \infty} J(u_n) \geq J(\bar{u}),
	\]
	which implies that $ J $ attained its infimum at $ \bar{u} $.
	
	It remains to show the uniqueness of $ \bar{u} $ under the condition that the function $ (u,\xi) \mapsto L(x,u,\xi) $ is strictly convex. Here we give the proof by contradiction.
	
	Assume that there exists $ \bar{v} $ different from $ \bar{u} $ but $ J(\bar{v}) = J(\bar{u}) = m $. By the strict convexity, for every $ t \in (0,1) $, we have that
	\begin{equation}\label{E4.14}
		t L(x,\bar{u},\nabla \bar{u}) + (1 - t)L(x,\bar{v},\nabla \bar{v}) - L\bigl(x,t \bar{u} + (1 - t) \bar{v},t \nabla \bar{u} + (1 - t) \nabla \bar{v}\bigr) > 0.
	\end{equation}
	Integrating both sides of the above inequality with respect to $ x $, we obtain
	\[
		t J(\bar{u}) + (1 - t) J(\bar{v}) - J\bigl(t \bar{u} + (1 - t) \bar{v}\bigr) \geq 0, \: \forall t \in (0,1).
	\]
	Since the integrand is non-negative, the inequality above becomes an equality if and only if the integrand is $ 0 $ almost everywhere, which is in contradiction to \eqref{E4.14}. Hence we have
	\[
		J(t \bar{u} + (1 - t) \bar{v}) < t J(\bar{u}) + (1 - t) J(\bar{v}) = m, \: \forall t \in (0,1).
	\]
	However, it contradicts the fact that $ m $ is the infimum. Thus, the uniqueness part has been proved. This completes the proof of the theorem.
\end{proof}

\begin{proof}[\textbf{Proof of \textup{Proposition \ref{Pr1}}.}]
	First, let us consider a simpler case. Suppose that there exist some $ a \in \mathbb{R}_{> 0} $ and $ b \in \mathbb{R}_{\geq 0} $ such that
	\[
		L(t,u,\dot{u}) = a \lvert \dot{u} \rvert + b \lvert u \rvert.
	\]
	We have
	\[
		J(u) = \int_0^1 L(t,u,\dot{u}) \mathrm{d}t = a \lVert \dot{u} \rVert_{L^1(0, 1)} + b \lVert u \rVert_{L^1(0, 1)}.
	\]
	
	Take $ u_0(t) = t $. Then we have that
	\[
		J(u_0) = a + \frac{b}{2} < + \infty.
	\]
	For any minimizing sequence $ \{u_n\} \subset E $, using Theorem \ref{T2}, $ \{u_n\} $ is bounded. Denote
    \[
    v_n(t) := \begin{cases}
        0, &t \in (0,1-\frac{1}{2^n}];\\
        u_n\bigl(1-2^n(1-x)\bigr), &t \in (1-\frac{1}{2^n},1).
    \end{cases}
    \]
    Obviously that $ \lVert \dot{v}_n \rVert_{L^1(0,1)} = \lVert \dot{u}_n \rVert_{L^1(0,1)} $ and $ \lVert v_n \rVert_{L^1(0,1)} = \frac{1}{2^n}\lVert u_n \rVert_{L^1(0,1)} $. 

    If $ b > 0 $, we have that $ a \leq J(v_n) < J(u_n) \to J(\bar{u}) $. Then $ \{v_n\} $ is also a minimizing sequence in $ E $. Since $ u_n $ is bounded, we obtain that $ \lVert v_n \rVert_{L^1(0,1)} \to 0 $ and
    \[
    \lim_{n \to +\infty} J(u_n) = \lim_{n \to +\infty} J(v_n) = \lim_{n \to +\infty} a\lVert \dot{u}_n \rVert_{L^1(0,1)}.
    \]
    And it implies that $ \lVert u_n \rVert_{L^1(0,1)} \to 0 $. For each $ n \in \mathbb{Z}_{>0}$, denote
    \[
    A_n := \{t\in (0,1) \vert \: u_n(t) \geq \frac{1}{2} \}.
    \]
    Since $ \lVert u_n \rVert_{L^1(0,1)} \to 0 $, the Lebesgue measure $ \lambda(A_n) \to 0 $. But $ \int_{A_n} \lvert \dot{u}_n \rvert \mathrm{d}t \geq \frac{1}{2} $, which contradicts Statement $ \mathbf{EI^{\lvert \alpha \rvert = 1,1}}{[(0,1), \{u_n\}, \exists \delta]} $ in Theorem \ref{T1'}. Hence the minimizing sequence $ \{u_n\} $ is not relatively weakly compact in $ W^{1,1}(0,1) $.

    If $ b = 0 $, using the Newton-Leibniz formula, there holds the inequality
    \[
    \lVert \dot{u} \rVert_{L^1(0, 1)} \geq \int_0^1 \dot{u}(t) \mathrm{d}t = u(1) - u(0) = 1,
    \]
		with equality if and only if $ \dot{u}(t) \geq 0 $ for almost everywhere in $ (0, 1) $. Obviously, every monotonically increasing function with respect to $ t $ minimizes the functional $ J $ already.
    
	In conclusion, the functional
	\[
		J(u) = a \lVert \dot{u} \rVert_{L^1(0, 1)} + b \lVert u \rVert_{L^1(0, 1)}
	\]
	has a minimal point $ \bar{u} $ if and only if $ b = 0 $. Moreover, we also obtain that $ J(\bar{u}) = a $ and $ \dot{u}(t) \geq 0 $ for almost everywhere in $ (0, 1) $.
	
	Next, we consider the case where the integrand has the form
	\[
		L(t,u,\dot{u}) = (k \lvert \dot{u} \rvert^p + \ell \lvert u \rvert^p)^{\frac{1}{p}}.
	\]
	Thus we have
	\[
		J(u) = \int_0^1 L(t,u,\dot{u}) \mathrm{d}t = \int_0^1 (k \lvert \dot{u} \rvert^p + \ell \lvert u \rvert^p)^{\frac{1}{p}} \mathrm{d}t.
	\]
	
	If $ 1 < p < + \infty $, let $ a = k^{\frac{1}{p}} > 0 $ and $ b = \ell^{\frac{1}{p}} \geq 0 $. 
	
	The following inequality will be used later:
	\begin{equation}\label{E4.36}
		x^p + y^p \leq (x + y)^p,
	\end{equation}
	with equality if and only if $ xy = 0 $, where $ x $, $ y \geq 0 $ and $ p > 1 $.

	We proceed to prove \eqref{E4.36}. Obviously, if $ xy = 0 $, the equal sign holds, so does \eqref{E4.36}. Besides, if $ xy \neq 0 $, denote $ k := \frac{x}{y} $. Notice that $ k > 0 $ and $ p > 1 $, thus
	\[
		\begin{split}
			\frac{x^p + y^p}{(x + y)^p} &= \frac{x^p}{(x + y)^p} + \frac{y^p}{(x + y)^p} = (\frac{\frac{x}{y}}{\frac{x}{y} + 1})^p + (\frac{1}{(\frac{x}{y} + 1)})^p\\
			&= (\frac{k}{k + 1})^p + (\frac{1}{k + 1})^p < \frac{k}{k + 1} + \frac{1}{k + 1} = 1.
		\end{split}
	\]
	Hence the proof of \eqref{E4.36} is completed.
	
	Take $ x = a \lvert \dot{u} \rvert $ and $ y = b \lvert u \rvert $ in \eqref{E4.36}. We have
	\begin{equation}\label{E4.43}
		a \lvert \dot{u} \rvert \leq (k \lvert \dot{u} \rvert^p + \ell \lvert u \rvert^p)^{\frac{1}{p}} = \bigl((a \lvert \dot{u} \rvert)^p + (b \lvert u \rvert)^p\bigr)^{\frac{1}{p}} \leq a \lvert \dot{u} \rvert + b \lvert u \rvert.
	\end{equation}

	Let
	\[
		I(u) = \int_0^1 a \lvert \dot{u} \rvert \mathrm{d}t \text{ and } H(u) = \int_0^1 (a \lvert \dot{u} \rvert + b \lvert u \rvert) \mathrm{d}t
	\]
	be two functionals defined on $ E $. According to the discussion above, the infimum of $ I $ is $ a $ and can be achieved. Since every sequence $ \{u_n\} $ minimizing $ H $ satisfies that $ \lVert u_n \rVert_{L^1(0,1)} \to 0 $, we have
	\[
		\inf_{u \in E} H(u) \leq \lim_{n \to +\infty} H(u_n) = a + \lim_{n \to +\infty} b\lVert u_n \rVert_{L^1(0,1)} = a.
	\]

	Integrating each part of \eqref{E4.43} with respect to $ t $, we have
	\begin{equation}\label{67}
		I(u) \leq J(u) \leq H(u), \: \forall u \in E.
	\end{equation}
	Hence
	\[
		a = \inf_{u \in E} I(u) \leq \inf_{u \in E} J(u) \leq \inf_{u \in E} H(u) \leq a.
	\]

	Let $ \{u_n\} $ be a minimizing sequence of $ J $ satisfying $ u_n \rightharpoonup \bar{u} $ and $ J(\bar{u}) = a $. According to \eqref{67}, $ I(\bar{u}) = a $. It follows that
	\[
		\int_0^1 (k \lvert \dot{\bar{u}} \rvert^p + \ell \lvert \bar{u} \rvert^p)^{\frac{1}{p}} \mathrm{d}t = a = \int_0^1 a \lvert \dot{\bar{u}} \rvert \mathrm{d}t = \int_0^1 (k \lvert \dot{\bar{u}} \rvert^p)^{\frac{1}{p}} \mathrm{d}t.
	\]
	Therefore, $ \ell = 0 $ is necessary for $ J $ to attain the infimum. It remains to show that $ \ell = 0 $ is sufficient too.

	If $ \ell = 0 $, we have
	\[
		J(u) = I(u) = \int_0^1 a \lvert \dot{u} \rvert \mathrm{d}t, \: \forall u \in E.
	\]
	Using the Newton-Leibniz formula, for every $ \bar{u} \in E $ satisfying that $ \dot{\bar{u}}(t) \geq 0 $ for almost everywhere in $ (0, 1) $, there holds $ J(\bar{u}) = a $. Thus, the sufficiency is proved.

	Finally, since $ \int_0^1 c(t) \mathrm{d}t $ is a constant, it does not affect whether the infimum of the functional $ J $ can be achieved. This completes the proof of Proposition \textup{\ref{Pr1}}.
\end{proof}

\backmatter

\bmhead{Acknowledgments} 

This work was supported partially by NSF of China, No. 12071316. We would also like to thank Professor Zhu Chaofeng from Nankai University for his discussion on whether interpolation inequalities in continuous differentiable function spaces can be applied. We also thank Professor Zhou Feng from East China Normal University for his fruitful comments about the boundary conditions of the domain $ \Omega $ in each theorem. M.J. would like to thank Professor Hongjie Dong from Brown University for helpful conversation.

\bmhead{Data availability statement}
This manuscript has no associated data.

\section*{Statements and Declarations}

\bmhead{Conflict of interest}
The authors declared that they have no conflicts of interest to this work.

\begin{appendices}

\section{Proofs of Corollaries}\label{Poc}

\begin{proof}[\textbf{Proof of \textup{Corollary~\ref{cor1}}.}]
	The proof of this corollary is almost ``isomorphic'' to the proof of Theorem \ref{T1}.\\

	\noindent Proof of ``$ \mathbf{H^{k, p}} \wedge \mathbf{B^{k, p}} \wedge \mathbf{EI_{\boldsymbol{\tau}}^{k, p}} \wedge \mathbf{E_{\infty}^{k, p}}[\mathbb{R}^d, \mathcal{F}, \exists M, \exists \delta, \exists \omega] \Rightarrow \mathbf{H^{k, p}} \wedge \mathbf{C^{k, p}}[\mathbb{R}^d, \mathcal{F}] $'': 

	Of course, we suppose that $ \mathbf{H^{k, p}} \wedge \mathbf{B^{k, p}} \wedge \mathbf{EI_{\boldsymbol{\tau}}^{k, p}} \wedge \mathbf{E_{\infty}^{k, p}}[\mathbb{R}^d, \mathcal{F}, \exists M, \exists \delta, \exists \omega] $. Then 
	\begin{itemize}
		\item[1.] For any $ \iota f \in \iota \mathcal{F} $, its norm satisfies
		\[
			\lVert \iota f \rVert_{L^p(\bigcup_{\lvert \alpha \rvert \leq k} \mathbb{R}^d_{\alpha})} = \sum_{\lvert \alpha \rvert \leq k}  \lVert D^{\alpha}f \rVert_{L^p(\mathbb{R}^d)} = \lVert f \rVert_{W^{k, p}(\mathbb{R}^d)} \leq M.
		\]
		\item[2.] For any $ \iota f \in \iota \mathcal{F} $, $ \forall \varepsilon \in \mathbb{R}_{> 0} $ and $ \forall h \in \mathbb{R}^d $ with $ \lVert h \rVert_{\mathbb{R}^d} < \delta(\varepsilon) $, we have
		\[
			\begin{split}
				\lVert \iota (\tau_h f) - \iota f \rVert_{L^p(\bigcup_{\lvert \alpha \rvert \leq k} \mathbb{R}^d_{\alpha})} &= \sum_{\lvert \alpha \rvert \leq k}  \lVert D^{\alpha}(\tau_h f) - D^{\alpha}f \rVert_{L^p(\mathbb{R}^d)}\\
				&= \lVert \tau_h f - f \rVert_{W^{k, p}(\mathbb{R}^d)} < \varepsilon.
			\end{split}
		\]
		\item[3.] Take $ \omega^{(k)} \colon \varepsilon \mapsto \bigcup_{\lvert \alpha \rvert \leq k} I_{\alpha} \circ \omega(\varepsilon) $. For any $ \iota f \in \iota \mathcal{F} $ and $ \forall \varepsilon \in \mathbb{R}_{> 0} $, we have
		\[
			\lVert \iota f \rVert_{L^p\bigl((\bigcup_{\lvert \alpha \rvert \leq k} \mathbb{R}^d_{\alpha}) \setminus \omega^{(k)}(\varepsilon)\bigr)} = (\sum_{\lvert \alpha \rvert \leq k}  \int_{\mathbb{R}^d \setminus \omega(\varepsilon)} \lvert D^{\alpha}f \rvert^p \mathrm{d}\lambda)^{\frac{1}{p}} < \varepsilon.
		\]
	\end{itemize}
	Using the Kolmogorov-M.Riesz-Fr\'echet theorem \cite[see][Corollary 4.27 on page 113]{Brezis}, we obtain that ``$ \iota \mathcal{F} $ is precompact in $ L^p(\bigcup_{\lvert \alpha \rvert \leq k} \mathbb{R}^d_{\alpha}) $'', which is equivalent to that ``$ \mathcal{F} $ is precompact in $ W^{k, p}(\mathbb{R}^d) $''.\\

	\noindent Proof of ``$ \mathbf{H^{k, p}} \wedge \mathbf{C^{k, p}}[\mathbb{R}^d, \mathcal{F}] \Rightarrow \mathbf{H^{k, p}} \wedge \mathbf{B^{k, p}} \wedge \mathbf{EI_{\boldsymbol{\tau}}^{k, p}} \wedge \mathbf{E_{\infty}^{k, p}}[\mathbb{R}^d, \mathcal{F}, \exists M, \exists \delta, \exists \omega] $'':

	Similarly, suppose we are given $ \mathbf{H^{k, p}} \wedge \mathbf{C^{k, p}}[\mathbb{R}^d, \mathcal{F}] $. By the definition of $ \iota $, we have that $ \iota \mathcal{F} $ is precompact in $ L^p(\bigcup_{\lvert \alpha \rvert \leq k} \mathbb{R}^d_{\alpha}) $. Using the Kolmogorov-M.Riesz-Fr\'echet theorem again, we obtain that
	\begin{itemize}
		\item[1.] There exists a positive constant $ M $ such that for any $ \iota f \in \iota \mathcal{F} $, its norm satisfies $ \lVert \iota f \rVert_{L^p(\bigcup_{\lvert \alpha \rvert \leq k} \mathbb{R}^d_{\alpha})} \leq M $. Thus we have
		\[
			\lVert f \rVert_{W^{k, p}(\mathbb{R}^d)} = \sum_{\lvert \alpha \rvert \leq k}  \lVert D^{\alpha}f \rVert_{L^p(\mathbb{R}^d)} = \lVert \iota f \rVert_{L^p(\bigcup_{\lvert \alpha \rvert \leq k} \mathbb{R}^d_{\alpha})} \leq M, 
		\]
		which implies $ \mathbf{B^{k, p}}[\mathbb{R}^d, \mathcal{F}, \exists M] $.
		\item[2.] There exists a function $ \delta \colon \mathbb{R}_{> 0} \to \mathbb{R}_{> 0} $ such that for any $ \iota f \in \iota \mathcal{F} $, $ \forall \varepsilon \in \mathbb{R}_{> 0} $ and $ \forall h \in \mathbb{R}^d $ with $ \lVert h \rVert_{\mathbb{R}^d} < \delta(\varepsilon) $, there holds $ \lVert \iota (\tau_h f) - \iota f \rVert_{L^p(\bigcup_{\lvert \alpha \rvert \leq k} \mathbb{R}^d_{\alpha})} < \varepsilon $. Then we have
		\[
			\begin{split}
				\lVert \tau_h f - f \rVert_{W^{k, p}(\mathbb{R}^d)} &= \sum_{\lvert \alpha \rvert \leq k}  \lVert D^{\alpha}(\tau_h f) - D^{\alpha}f \rVert_{L^p(\mathbb{R}^d)}\\
				&= \lVert \iota (\tau_h f) - \iota f \rVert_{L^p(\bigcup_{\lvert \alpha \rvert \leq k} \mathbb{R}^d_{\alpha})} < \varepsilon.
			\end{split}
		\]
		Thus we get $ \mathbf{EI_{\boldsymbol{\tau}}^{k, p}}[\mathbb{R}^d, \mathcal{F}, \exists \delta] $.
		\item[3.] There exists a set-valued mapping $ \omega^{(k)} \colon \mathbb{R}_{> 0} \to 2^{\bigcup_{\lvert \alpha \rvert \leq k} \mathbb{R}^d_{\alpha}} $ such that for every $ \iota f \in \iota \mathcal{F} $ and $ \forall \varepsilon \in \mathbb{R}_{> 0} $, there holds $ \lVert \iota f \rVert_{L^p\bigl((\bigcup_{\lvert \alpha \rvert \leq k} \mathbb{R}^d_{\alpha}) \setminus \omega^{(k)}(\varepsilon)\bigr)} < \varepsilon $. Take $ \omega \colon \varepsilon \mapsto \bigcup_{\lvert \alpha \rvert \leq k} I_{\alpha}^{- 1}\bigl(\omega^{(k)}(\varepsilon) \cap \mathbb{R}_{\alpha}\bigr) $, then we have
		\[
			\begin{split}
				(\sum_{\lvert \alpha \rvert \leq k}  \int_{\mathbb{R}^d \setminus \omega(\varepsilon)} \lvert D^{\alpha}f \rvert^p \mathrm{d}\lambda)^{\frac{1}{p}} &= (\int_{\bigcup_{\lvert \alpha \rvert \leq k} \Bigl(\mathbb{R}^d_{\alpha} \setminus \bigl(I_{\alpha} \circ \omega(\varepsilon)\bigr)\Bigr)} \lvert \iota f \rvert^p \mathrm{d}\mu)^{\frac{1}{p}}\\
				&\leq \lVert \iota f \rVert_{L^p\bigl((\bigcup_{\lvert \alpha \rvert \leq k} \mathbb{R}^d_{\alpha}) \setminus \omega^{(k)}(\varepsilon)\bigr)} < \varepsilon,
			\end{split}
		\]
		which leads to $ \mathbf{E_{\infty}^{k, p}}[\mathbb{R}^d, \mathcal{F}, \exists \omega] $ immediately.
	\end{itemize}

	The proof is completed.
\end{proof}

\begin{proof}[\textbf{Proof of \textup{Corollary~\ref{cor1'}}.}]
    The proof of this corollary is almost ``isomorphic'' to the proof of Theorem \ref{T1'}. Similarly, we only need to prove the following two.
    
    \begin{description}
        \item[\textup{1.}] $ \mathbf{H^{k, p}_{b.o.0}} \wedge \mathbf{B^{\lvert \alpha \rvert \in \{0, k\},p}}{[\Omega, \mathcal{F}, \exists M]} \Rightarrow \mathbf{H^{k, p}_{b.o.0}} \wedge \mathbf{B^{k,p}}{[\Omega, \mathcal{F}, \exists M']} $. 
    \end{description}
    
    It can be obtained directly from the interpolation inequality \cite[see][Theorem 5.2 on page 135]{Adams}.
    
    \begin{description}  
        \item[\textup{2.}] $ \mathbf{H^{k, p}_{b.o.0}}\wedge \mathbf{B^{\lvert \alpha \rvert \in \{0, k\},p}} \wedge \mathbf{EI_{\boldsymbol{\tau}}^{\lvert \alpha \rvert = k,p}}{[\Omega, \mathcal{F}, \exists M, \exists \delta]} \Rightarrow \mathbf{H^{k, p}_{b.o.0}}\wedge \mathbf{B^{\lvert \alpha \rvert \in \{0, k\},p}} \wedge \mathbf{EI_{\boldsymbol{\tau}}^{k,p}}{[\Omega, \mathcal{F}, \exists M', \exists \delta']} $.
    \end{description}

        From the above discussion we already know that $ \mathcal{F} $ is bounded in $ W^{k,p}(\Omega) $. Applying the Rellich-Kondrachov theorem \cite[see][Theorem 6.3 on page 168]{Adams}, we obtain that the embedding from $ W_0^{k,p}(\Omega) $ to $ W^{k-1,p}(\Omega) $ is compact (When $ k = 1 $, the notation $ W^{0,p} $ is regarded as $ L^p $). Hence $ \iota\mathcal{F} $ is precompact in $ L^p(\Omega^{(k)}) $. Then applying Corollary \ref{cor1}, the statement $ \mathbf{EI_{\boldsymbol{\tau}}^{k,p}}{[\Omega, \mathcal{F}, \exists M', \exists \delta']} $ is satisfied.
        
    Thus the proof is completed.
\end{proof}

\begin{proof}[\textbf{Proof of \textup{Corollary~\ref{cor2}}.}]
	We need to fine-tune the definition of the isometric operator $ \iota $ so that it "translates" from the Sobolev space to the space of continuously differentiable functions. Denote 
	\[
		K^{(m)} := \bigsqcup_{\lvert \alpha \rvert \leq m} K = \bigcup_{\lvert \alpha \rvert \leq m} K_{\alpha},
	\]
	where $ K_{\alpha} := I_{\alpha}(K) $ is a disjoint copy of $ K $ given by the isometric mapping $ I_{\alpha} $. Moreover, suppose $ K^{(m)} $ satisfies that for each fixed $ \alpha $ there holds
	\[
		\mathrm{diam}(K) := \max_{(x_1, x_2) \in K \times K} d(x_1, x_2) < \min_{(x_{\alpha}, x_{\not\alpha}) \in K_{\alpha} \times (K^{(m)} \setminus K_{\alpha})} d(x_{\alpha}, x_{\not\alpha}) < + \infty,
	\]
	i.e., the distance between two different ``$ K_{\alpha} $''s is greater than the diameter of $ K $ but less than $ + \infty $, which is always easy to satisfy. By these assumptions, we have that $ K^{(m)} $ is a compact metric space. Then we denote
	\[
		\begin{split}
			\iota \colon \quad C^m(K) &\to C(K^{(m)})\\
			u &\mapsto \iota u,
		\end{split}
	\]
	where
	\[
		\begin{split}
			\iota u \colon \quad K^{(m)} &\to \mathbb{R}\\
			y &\mapsto (\bigoplus_{\lvert \alpha \rvert \leq m} (D^{\alpha}u) \circ I_{\alpha}^{- 1})(y)\\
			i.e. \: y &\mapsto (D^{\alpha}u) \circ I_{\alpha}^{- 1}(y), \: \text{if } y \in K_{\alpha}.
		\end{split}
	\]
	It is easy to verify that $ \iota $ is well-defined. Let $ W $ be the range of $ \iota $, so $ W $ is a closed subspace of $ C(K^{(m)}) $. The same arguments as in Lemma \ref{lemmaIs} remind us $ \iota $ is an isometry from $ C^m(K) $ to $ W $. We denote its inverse operator as $ \iota ^{- 1} \colon W \to C^m(K) $.\\

	\noindent Proof of ``$ \mathbf{H^m} \wedge \mathbf{B^m} \wedge \mathbf{EC^m}[K, \mathcal{F}, \exists M, \exists \delta] \Rightarrow \mathbf{H^m} \wedge \mathbf{C^m}[K, \mathcal{F}] $'': 

	Suppose $ \mathbf{H^m} \wedge \mathbf{B^m} \wedge \mathbf{EC^m}[K, \mathcal{F}, \exists M, \exists \delta] $. Then 
	\begin{itemize}
		\item[1.] For any $ \iota f \in \iota \mathcal{F} $, its norm satisfies
		\[
			\lVert \iota f \rVert_{C(K^{(m)})} = \max_{\lvert \alpha \rvert \leq m} \lVert D^{\alpha}f \rVert_{C(K)} = \lVert f \rVert_{C^m(K)} \leq M.
		\]
		\item[2.] Take $ \delta^{(m)} \colon \varepsilon \mapsto \min_{\lvert \alpha \rvert \leq k} \{\delta(\varepsilon), \mathrm{diam}(K)\} $. For any $ \iota f \in \iota \mathcal{F} $, $ \forall \varepsilon \in \mathbb{R}_{> 0} $ and $ \forall y_1,y_2 \in K^{(m)} $ with $ d(y_1, y_2) < \delta^{(m)}(\varepsilon) $, since $ y_1 $ and $ y_2 $ must come from a same ``$ K_{\alpha} $'' (without loss of generality, we name it $ K_{\alpha'} $), thus we have
		\[
			\lvert \iota f(y_1) - \iota f(y_2) \rvert \leq \max_{\lvert \alpha \rvert \leq m} \lvert (D^{\alpha}f)\bigl(I_{\alpha'}^{- 1}(y_1)\bigr) - (D^{\alpha}f)\bigl(I_{\alpha'}^{- 1}(y_2)\bigr) \rvert < \varepsilon.
		\]
	\end{itemize}
	Using the Ascoli-Arzel\`a theorem \cite[see][Theorem 4.25 on page 111]{Brezis}, we obtain that ``$ \iota \mathcal{F} $ is precompact in $ C(K^{(m)}) $'', which is equivalent to that ``$ \mathcal{F} $ is precompact in $ C^m(K) $''.\\

	\noindent Proof of ``$ \mathbf{H^m} \wedge \mathbf{C^m}[K, \mathcal{F}] \Rightarrow \mathbf{H^m} \wedge \mathbf{B^m} \wedge \mathbf{EC^m}[K, \mathcal{F}, \exists M, \exists \delta] $'': 

	Suppose that we are given $ \mathbf{H^m} \wedge \mathbf{C^m}[K, \mathcal{F}] $. By the definition of $ \iota $, we have that $ \iota \mathcal{F} $ is precompact in $ C(K^{(m)}) $. Using the Ascoli-Arzel\`a theorem again, we obtain that
	\begin{itemize}
		\item[1.] There exists a positive constant $ M $ such that for any $ \iota f \in \iota \mathcal{F} $, its norm satisfies $ \lVert \iota f \rVert_{C(K^{(m)})} \leq M $. Thus we have
		\[
			\lVert f \rVert_{C^m(K)} = \max_{\lvert \alpha \rvert \leq m} \lVert D^{\alpha}f \rVert_{C(K)} = \lVert \iota f \rVert_{C(K^{(m)})} \leq M, 
		\]
		which implies $ \mathbf{B^m}[K, \mathcal{F}, \exists M] $.
		\item[2.] There exists a function $ \delta \colon \mathbb{R}_{> 0} \to \mathbb{R}_{> 0} $ such that for any $ \iota f \in \iota \mathcal{F} $, $ \forall \varepsilon \in \mathbb{R}_{> 0} $ and $ \forall x_1, \: x_2 \in K $ with $ d(x_1, x_2) < \delta(\varepsilon) $ (note that $ d\bigl(I_{\alpha}(x_1), I_{\alpha}(x_2)\bigr) = d(x_1, x_2) < \delta(\varepsilon) $), there holds $ \lvert \iota f \circ I_{\alpha}(x_1) - \iota f \circ I_{\alpha}(x_2) \rvert < \varepsilon $ for every $ \alpha $. Then we have
		\[
			\max_{\lvert \alpha \rvert \leq m} \lvert (D^{\alpha}f)(x_1) - (D^{\alpha}f)(x_2) \rvert = \max_{\lvert \alpha \rvert \leq m} \lvert \iota f \circ I_{\alpha}(x_1) - \iota f \circ I_{\alpha}(x_2) \rvert < \varepsilon.
		\]
		Thus we get $ \mathbf{EC^m}[K, \mathcal{F}, \exists \delta] $.
	\end{itemize}

	The proof is completed.
\end{proof}

\begin{proof}[\textbf{Proof of \textup{Corollary~\ref{cor2'}}.}]
    We only need to show
    \begin{description}
        \item $ \mathbf{H^m} \wedge \mathbf{B^m} \wedge \mathbf{EC^{\lvert \alpha \rvert = m}}{[K, \mathcal{F}, \exists M, \exists \delta]} \Rightarrow \mathbf{H^m} \wedge \mathbf{B^m} \wedge \mathbf{EC^m}[K, \mathcal{F}, \exists M, \exists \delta'] $.
    \end{description}

    Note that the embedding from $ C^1(K) $ to $ C(K) $ is compact. It follows that the embedding from $ C^m(K) $ to $ C^{m-1}(K) $ is compact. Hence $ \iota \mathcal{F} $ is precompact in $ C(K^{(m)}) $. Then applying Corollary \ref{cor2}, the statement $ \mathbf{EC^m}[K, \mathcal{F}, \exists M, \exists \delta'] $ is satisfied.

    The proof is completed.
\end{proof}

\section{Syntactic Category}\label{Details}
To describe and prove Theorem \ref{T1} more clearly, we introduce a few notations from category theory, but we will not be using some overly abstract concepts. Inspired by the syntactic category \cite[see][page 10]{John}, we simplify sentences into symbols and make the following definition, which is different from the original meaning.

\begin{definition}\label{category}
	Let $ S $ be a set of sentences in some propositional language. For simplicity in this paper, we encourage the reader to think of $ S $ as ``a collection of mathematical statements with some axioms". We define the \textit{synatic category determined by $S$} as the category $\mathfrak{C}(S)$ with the following data:
	\begin{description}
		\item[\textbf{\textup{Objects:}}] 
		The sentences in $ S $ and those formed by connecting some elements in $ S $ with the connectives $ \wedge $ any number of times. In particular, in this paper, all sentences represented by abbreviations are elements of $S$. For example, $ \mathbf{B^1}{[\Omega, \mathcal{F}, \exists M]} $, $ \mathbf{EI^1}{[\Omega, \mathcal{F}, \exists \delta]} $, $ \mathbf{E_{\infty}^1}{[\Omega, \mathcal{F}, \exists \omega]} $ and $ \mathbf{B^1} \wedge \mathbf{EI^1} \wedge \mathbf{E_{\infty}^1}{[\Omega, \mathcal{F}, \exists M, \exists \delta, \exists \omega]} $ (in Theorem \ref{Dunford-Pettis}) are all in $S$.

		\item[\textbf{\textup{Morphisms:}}] 
		For objects $ p $ and $ q $, the morphism $ p \overset{P}{\Rightarrow} q $ coincides with a proof $ P $ from $ p $ to $ q $. For example, the Bolzano-Weierstrass (BW) theorem gives a morphism ``$\overset{\text{BW}}{\Rightarrow}$'' from the sentence ``$\text{Set } A \text{ is bounded and closed in } \mathbb{R} $'' to the sentence ``$ A \text{ is a compact set} $''.\\
  
        We assume that each sentence can be proved from itself without any reason, that is, $ p \overset{\varnothing}{\Rightarrow} p $ for every $ p $. Since we only focus on the existence of proofs, we will omit the information about what the proof is actually on the arrows, that is, we treat different proofs e.g. $ p \overset{P}{\Rightarrow} q $ and $ p \overset{Q}{\Rightarrow} q $ as a same morphism $ p \Rightarrow q $ hereafter. In other words, $ \hom(p, q) $ contains at most one element for every $ p $ and $ q $. It also means that $ \mathfrak{C}(S) $ is a thin category.
		\item[\textbf{\textup{Composites:}}] 
		given two morphisms $ p \Rightarrow q $ and $ q \Rightarrow r $, the composite $ p \Rightarrow r $ is naturally the proof obtained by concatenating the two proofs. %For example: ``$ A $ is a compact set'' $ \wedge $ ``$ f \colon A \to \mathbb{R} $ is a continuous function'' $ \Rightarrow $ ``$ f(A) $ is a compact set in $ \mathbb{R} $''$ \Rightarrow$``$ f(A) $ is bounded and closed in $ \mathbb{R} $''.
	\end{description}

Note that unlike the original definition \cite[see][page 10]{John}, we do not require that the sentences of the proofs must be in $ S $, that is, we can use any mathematical theorem in a self-consistent axiom system as a morphism. Therefore, what elements the set $S$ contains does not affect whether there exists a morphism between two objects in $S$. The following categorical data may be interpreted in $\mathfrak{C}$ as follows:
	\begin{description}
		\item[\textbf{\textup{Isomorphisms:}}] 
		an isomorphism is a morphism $ p \Rightarrow p' $ such that there exists a morphism $ p' \Rightarrow p $ as the inverse. We denote isomorphisms by $ p \Leftrightarrow p' $.
		\item[\textbf{\textup{General products:}}] \cite[see][page 22]{Abramsky} 
		for a family of objects $ \{p_{\alpha}\}_{\alpha \in \Lambda} $, a product is an object $ \bigwedge_{\alpha \in \Lambda} p_{\alpha} $ with canonical proofs of conjuncts from conjunctions as projection morphisms $ \bigwedge_{\alpha \in \Lambda} p_{\alpha} \Rightarrow p_{\alpha} $ for each $ p_{\alpha} $. Products exist in $ \mathfrak{C}(S) $ universally.
		\item[\textbf{\textup{Uniqueness up to unique isomorphism:}}] \cite[see][page 23]{Abramsky} 
		suppose that $ p $ and $ p' $ are both general products of objects $ \{p_{\alpha}\}_{\alpha \in \Lambda} $, then $ p $ and $ p' $ are isomorphic, i.e., $ p \Leftrightarrow p' $.
	\end{description}
\end{definition}

\begin{remark}
	Dually, one may define coproducts using the connectives $ \vee $ (i.e. or).
\end{remark}

\begin{remark}
	When we add a new sentence $ \gamma $ into $ S $ such that the category $ \mathfrak{C}(S) $ expands to $ \mathfrak{C}(S \cup \{\gamma\}) $, one may observe that $ \mathfrak{C}(S) $ is a full subcategory of $ \mathfrak{C}(S \cup \{\gamma\}) $, therefore morphisms in $ \mathfrak{C}(S) $ are inherited into $ \mathfrak{C}(S \cup \{\gamma\}) $. We can always make $ S $ contain enough sentences, so we do not strictly define what $ S $ contains hereafter.
\end{remark}

\newpage
\section{Notations and Abbreviations}\label{notations}

The following notations and abbreviations are used in this manuscript:

\begingroup
\renewcommand{\arraystretch}{1.5}
\noindent 
\begin{tabular}{p{0.23\linewidth} | p{0.72\linewidth}}
$\lambda$ & The Lebesgue measure in real Euclidean spaces\\
$\wedge$ & And\\
$ \Leftrightarrow $ & Is equivalent to\\
$ \Rightarrow $ & Deduces\\
$\alpha$ & The multi-index (see Notation \ref{multi-indices})\\
$\lvert \alpha \rvert$ & The length of multi-index $\alpha$\\
$ D^{\alpha}f $ & The weak generalized partial derivative of $ f $\\
%$ \tau_h f(x) $ & $f(x+h)$\\
$\Omega_\alpha$ & A copy of $\Omega$ with respect to $\alpha$ (see Definition \ref{definitionOmegak})\\
$I_\alpha$ & An isometry from $\Omega$ to $\Omega_\alpha$\\
$\lambda_\alpha$ & The Lebesgue measure in $\Omega_\alpha$\\
$\Omega^{(k)}$ & The disjoint union of all the $\Omega_\alpha$ with $\lvert \alpha \rvert \leq k$\\
$\mu$ & The measure in $\Omega^{(k)}$\\
$ \iota $ & A special isometry (see Definition \ref{definitionIs})
\end{tabular}
\vspace{0.5cm}\\
\noindent 
\begin{tabular}{p{0.23\linewidth} | p{0.72\linewidth}}
$ \mathbf{H} $ & Sentences for hypotheses\\
$ \mathbf{H_{\boldsymbol{\sigma}}^1}{[\Omega, \mathcal{F}]} $ & $ \Omega $ is a $ \sigma $-finite measure space, and $ \mathcal{F} $ is a subset of $ L^1(\Omega) $.\\
$ \mathbf{H^{k, 1}_{n.o.}}{[\Omega, \mathcal{F}]} $ & $ \Omega $ is a non-empty open set in $ \mathbb{R}^d $, and $ \mathcal{F} $ is a subset of $ W^{k,1}(\Omega) $.\\
$ \mathbf{H^{k, 1}_{b.o.cone}}{[\Omega, \mathcal{F}]} $ & $ \Omega $ is a bounded open set in $ \mathbb{R}^d $, and $ \Omega $ satisfies the cone condition. $ \mathcal{F} $ is a subset of $ W^{k,1}(\Omega) $.\\
$ \mathbf{H^{k, p}}{[\mathbb{R}^d, \mathcal{F}]} $ & $ \mathcal{F} $ is a subset of $ W^{k,p}(\mathbb{R}^d) $ with $ 1 \leq p < + \infty $.\\
$ \mathbf{H^m}{[K, \mathcal{F}]} $ & $K$ is a compact metric space which satisfies the uniform $ C^m $ regularity condition, and $ \mathcal{F} $ is a subset of $ C^m(K) $.
\end{tabular}

\noindent 
\begin{tabular}{p{0.30\linewidth} | p{0.65\linewidth}}
$ \mathbf{B} $ & Statement related to boundedness\\
$ \mathbf{B^1}{[\Omega, \mathcal{F}, \exists M]} $ & There exists $ M \in \mathbb{R}_{> 0} $ such that $ \lVert f \rVert_{L^1(\Omega)} \leq M $ for every $ f \in \mathcal{F} $.\\
$ \mathbf{B^{k, p}}{[\Omega, \mathcal{F}, \exists M]} $ & There exists $ M \in \mathbb{R}_{> 0} $ such that $ \lVert f \rVert_{W^{k,p}(\Omega)} \leq M $ for every $ f \in \mathcal{F} $.\\
$ \mathbf{B^{\lvert \alpha \rvert \in \{0, k\},p}}$ ${[\Omega, \mathcal{F}, \exists M]} $ & There exists a constant $ M \in \mathbb{R}_{> 0} $ such that 
                \[
				\lVert f \rVert_{L^p(\Omega)} + \sum_{\lvert \alpha \rvert = k} \lVert D^{\alpha} f \rVert_{L^p(\Omega)} \leq M
			\]
                for every $ f \in \mathcal{F} $.\\
$ \mathbf{B^m}{[K, \mathcal{F}, \exists M]} $ & There exists a constant $ M \in \mathbb{R}_{> 0} $ such that $ \max_{\lvert \alpha \rvert \leq m} \lVert D^{\alpha}f \rVert_{C(K)} \leq M $ for every $ f \in \mathcal{F} $.
\end{tabular}
\vspace{0.5cm}\\
\noindent 
\begin{tabular}{p{0.30\linewidth} | p{0.65\linewidth}}
$ \mathbf{EC} $ & Statement related to equi-continuity.\\
$ \mathbf{EC^m}{[K, \mathcal{F}, \exists \delta]} $ & There exists a function $ \delta \colon \mathbb{R}_{> 0} \to \mathbb{R}_{> 0} $ such that
			\[
				\max_{\lvert \alpha \rvert \leq m} \lvert (D^{\alpha}f)(x_1) - (D^{\alpha}f)(x_2) \rvert < \varepsilon
			\]
			for every $ f \in \mathcal{F} $, for any $ \varepsilon \in \mathbb{R}_{> 0} $, and for all $ x_1 $ and $ x_2 $ in $ K $ with the distance between them $ d(x_1, x_2) < \delta(\varepsilon) $.\\
$ \mathbf{EC^{\lvert \alpha \rvert = m}}{[K, \mathcal{F}, \exists \delta]} $ & There exists a function $ \delta \colon \mathbb{R}_{> 0} \to \mathbb{R}_{> 0} $ such that
			\[
				\max_{\lvert \alpha \rvert = m} \lvert (D^{\alpha}f)(x_1) - (D^{\alpha}f)(x_2) \rvert < \varepsilon
			\]
			for every $ f \in \mathcal{F} $, for any $ \varepsilon \in \mathbb{R}_{> 0} $, and for all $ x_1 $ and $ x_2 $ in $ K $ with the distance between them $ d(x_1, x_2) < \delta(\varepsilon) $.
\end{tabular}

\noindent 
\begin{tabular}{p{0.23\linewidth} | p{0.72\linewidth}}
$ \mathbf{EI} $ & Statement related to equi-integrability.\\
$ \mathbf{EI^1}{[\Omega, \mathcal{F}, \exists \delta]} $. & There exists a function $\delta \colon \mathbb{R}_{> 0} \to \mathbb{R}_{> 0}$ such that $\int_A \lvert f \rvert \mathrm{d}\lambda < \varepsilon$ for every $f \in \mathcal{F}$, for any $\varepsilon \in \mathbb{R}_{> 0}$, and for all measurable set $A \subseteq \Omega$ with its measure $\lambda(A) < \delta(\varepsilon)$.\\
$ \mathbf{EI^{k, 1}}{[\Omega, \mathcal{F}, \exists \delta]} $ & There exists a function $ \delta \colon \mathbb{R}_{> 0} \to \mathbb{R}_{> 0} $ such that
			\[
				\sum_{\lvert \alpha \rvert \leq k} \int_A \lvert D^{\alpha} f \rvert \mathrm{d}\lambda < \varepsilon
			\]
			for every $ f \in \mathcal{F} $, for any $ \varepsilon \in \mathbb{R}_{> 0} $, and for all $ A $ which is measurable in $ \Omega $ with its measure $ \lambda(A) < \delta(\varepsilon) $.\\
$ \mathbf{EI^{\lvert \alpha \rvert = k,1}}{[\Omega, \mathcal{F}, \exists \delta]} $ & There exists a function $ \delta \colon \mathbb{R}_{> 0} \to \mathbb{R}_{> 0} $ such that
			\[
				\sum_{\lvert \alpha \rvert = k} \int_A \lvert D^{\alpha} f \rvert \mathrm{d}\lambda < \varepsilon
			\]
			for every $ f \in \mathcal{F} $, for any $ \varepsilon \in \mathbb{R}_{> 0} $, and for all $ A $ which is measurable in $ \Omega $ with its measure $ \lambda(A) < \delta(\varepsilon) $.\\
$ \mathbf{EI_{\boldsymbol{\tau}}^{k, p}}{[\mathbb{R}^d, \mathcal{F}, \exists \delta]} $ & There exists a function $ \delta \colon \mathbb{R}_{> 0} \to \mathbb{R}_{> 0} $ such that
			\[
				\lVert \tau_h f - f \rVert_{W^{k, p}(\mathbb{R}^d)} < \varepsilon
			\]
			for every $ f \in \mathcal{F} $, for any $ \varepsilon \in \mathbb{R}_{> 0} $, and for all $ h \in \mathbb{R}^d $ with its norm $ \lVert h \rVert_{\mathbb{R}^d} < \delta(\varepsilon) $.\\
$ \mathbf{EI_{\boldsymbol{\tau}}^{\lvert \alpha \rvert = k,p}}{[\Omega, \mathcal{F}, \exists \delta]} $ & There exists a function $ \delta \colon \mathbb{R}_{> 0} \to \mathbb{R}_{> 0} $ such that
			\[
				\sum_{\lvert \alpha \rvert = k} \lVert \tau_h (D^{\alpha}f) - D^{\alpha}f \rVert_{L^p(\mathbb{R}^d)} < \varepsilon
			\]
			for every $ f \in \mathcal{F} $, for any $ \varepsilon \in \mathbb{R}_{> 0} $, and for all $ h \in \mathbb{R}^d $ with its norm $ \lVert h \rVert_{\mathbb{R}^d} < \delta(\varepsilon) $. Here, the functions $ D^{\alpha}f $ are extended to be $0$ outside $\Omega$.
\end{tabular}

\noindent 
\begin{tabular}{p{0.23\linewidth} | p{0.72\linewidth}}
$ \mathbf{E_{\infty}} $ & Sentences for some equi-properties at infinity\\
$\mathbf{E_{\infty}^{k,p}}{[\Omega, \mathcal{F}, \exists \omega]}$ & There exists a set-valued mapping $ \omega \colon \mathbb{R}_{> 0} \to 2^{\mathbb{R}^d} $ such that $ \omega(\varepsilon) $ is measurable with its measure $ \lambda \circ \omega(\varepsilon) $ being finite, and there satisfies that
			\[
				(\sum_{\lvert \alpha \rvert \leq k} \int_{\mathbb{R}^d \setminus \omega(\varepsilon)} \lvert D^{\alpha} f \rvert^p \mathrm{d}\lambda)^{\frac{1}{p}} < \varepsilon
			\]
	for every $ f \in \mathcal{F} $ and any $ \varepsilon \in \mathbb{R}_{> 0} $.\\
$ \mathbf{E_{\infty}^{k, p}}{[\Omega, \mathcal{F}, \exists \omega]} $ & There exists a set-valued mapping $ \omega \colon \mathbb{R}_{> 0} \to 2^{\Omega} $ such that $ \omega(\varepsilon) $ is measurable with its measure $ \lambda \circ \omega(\varepsilon) $ being finite, and there satisfies that
			\[
				(\sum_{\lvert \alpha \rvert \leq k} \int_{\Omega \setminus \omega(\varepsilon)} \lvert D^{\alpha} f \rvert^p \mathrm{d}\lambda)^{\frac{1}{p}} < \varepsilon
			\]
			for every $ f \in \mathcal{F} $ and any $ \varepsilon \in \mathbb{R}_{> 0} $.
\end{tabular}
\vspace{0.5cm}\\
\noindent 
\begin{tabular}{p{0.23\linewidth} | p{0.72\linewidth}}
$ \mathbf{C} $ & Statement related to compactness\\
$ \mathbf{C_w^1}{[\Omega, \mathcal{F}]} $ & $ \mathcal{F} $ is relatively weakly compact in $ L^1(\Omega) $.\\
$ \mathbf{C_w^{k, 1}}{[\Omega, \mathcal{F}]} $ & $ \mathcal{F} $ is relatively weakly compact in $ W^{k,1}(\Omega) $.\\
$ \mathbf{C^{k, p}}{[\mathbb{R}^d, \mathcal{F}]} $ & $ \mathcal{F} $ is precompact in $ W^{k, p}(\mathbb{R}^d) $.\\
$ \mathbf{C^m}{[K, \mathcal{F}]} $ & $ \mathcal{F} $ is precompact in $ C^m(K) $.
\end{tabular}

\endgroup

\begin{notation}[multi-indices $ \alpha $]\textup{\cite[3.5 on page 61]{Adams}}\label{multi-indices}
		Given integers $ d \geq 1 $ and $ k \geq 0 $, let
		\[
			\alpha := (\alpha_1, \alpha_2, \ldots, \alpha_d) \in \mathbb{N}^d
		\]
		be a multi-index such that
		\[
			\lvert \alpha \rvert := \sum_{i = 1}^d \alpha_i \leq k.
		\]
		For any $ f \in W^{k, p}(\Omega) $ where $ \Omega \subset \mathbb{R}^d $, we write
		\[
			D^{\alpha}f = \frac{\partial^{\alpha_1}}{\partial x_1^{\alpha_1}} \frac{\partial^{\alpha_2}}{\partial x_2^{\alpha_2}} \cdots \frac{\partial^{\alpha_d}}{\partial x_d^{\alpha_d}} f
		\]
		to denote the weak generalized partial derivative of $ f $ with respect to $\alpha$.
	\end{notation}

%An appendix contains supplementary information that is not an essential part of the text itself but which may be helpful in providing a more comprehensive understanding of the research problem or it is information that is too cumbersome to be included in the body of the paper.

%%=============================================%%
%% For submissions to Nature Portfolio Journals %%
%% please use the heading ``Extended Data''.   %%
%%=============================================%%

%%=============================================================%%
%% Sample for another appendix section			       %%
%%=============================================================%%

%% \section{Example of another appendix section}\label{secA2}%
%% Appendices may be used for helpful, supporting or essential material that would otherwise 
%% clutter, break up or be distracting to the text. Appendices can consist of sections, figures, 
%% tables and equations etc.

\end{appendices}

%%===========================================================================================%%
%% If you are submitting to one of the Nature Portfolio journals, using the eJP submission   %%
%% system, please include the references within the manuscript file itself. You may do this  %%
%% by copying the reference list from your .bbl file, paste it into the main manuscript .tex %%
%% file, and delete the associated \verb+\bibliography+ commands.                            %%
%%===========================================================================================%%

\bibliography{sn-bibliography}

\end{document}